\documentclass[11pt]{amsart}

\usepackage{amsfonts,amsmath,amssymb}

\newtheorem{theorem}{Theorem}[section]
\newtheorem{proposition}[theorem]{Proposition}
\newtheorem{lemma}[theorem]{Lemma}
\newtheorem{conjecture}[theorem]{Conjecture}
\newtheorem{corollary}[theorem]{Corollary}

\numberwithin{equation}{section}

\newcommand{\qbinsb}[2]{\biggl[\genfrac{}{}{0pt}{1}{#1}{#2}\biggr]}
\newcommand{\qbin}[2]{\genfrac{[}{]}{0pt}{}{#1}{#2}}
\newcommand{\qbins}[2]{{\genfrac{[}{]}{0pt}{1}{#1}{#2}}}
\newcommand{\ip}[2]{(#1|#2)}
\newcommand{\Z}{\mathbb{Z}}
\newcommand{\Pa}{\mathcal{P}}
\newcommand{\A}{\mathcal{A}}
\newcommand{\C}{\mathcal{C}}

\newcommand{\Su}{\mathcal{S}}
\newcommand{\pmd}[1]{\negthickspace\negthickspace\pod{#1}}

\begin{document}

\title{The Bailey lemma and Kostka polynomials}

\author[Ole Warnaar]{S. Ole Warnaar}

\address{Department of Mathematics and Statistics,
The University of Melbourne, Vic 3010, Australia}
\email{warnaar@ms.unimelb.edu.au}


\thanks{Work supported by the Australian Research Council}

\begin{abstract}
Using the theory of Kostka polynomials, we prove an A$_{n-1}$
version of Bailey's lemma at integral level.
Exploiting a new, conjectural expansion for Kostka numbers,
this is then generalized to fractional levels, leading to a new 
expression for admissible characters of A$_{n-1}^{(1)}$ and to 
identities for A-type branching functions.
\end{abstract}

\maketitle

\section{Introduction}
Let $\alpha=(\alpha_L)_{L\geq 0}$ and $\beta=(\beta_L)_{L\geq 0}$
be sequences and $(a;q)_r=(a)_r=\prod_{k=0}^{r-1}(1-aq^k)$ a $q$-shifted 
factorial. Then $(\alpha,\beta)$ is called a Bailey pair relative to $a$ if
\begin{equation}\label{BP}
\beta_L=\sum_{r=0}^L \frac{\alpha_r}{(q)_{L-r}(aq)_{L+r}}.
\end{equation}
Similarly, a pair of sequences $(\gamma,\delta)$ is called a conjugate 
Bailey pair relative to $a$ if
\begin{equation}\label{CBP}
\gamma_L=\sum_{r=L}^{\infty} \frac{\delta_r}{(q)_{r-L}(aq)_{r+L}}.
\end{equation}
By a simple interchange of sums it follows that if $(\alpha,\beta)$ and
$(\gamma,\delta)$ are a Bailey pair and conjugate Bailey pair relative to
$a$, then
\begin{equation}\label{abcd}
\sum_{L=0}^{\infty}\alpha_L\gamma_L=
\sum_{L=0}^{\infty}\beta_L\delta_L,
\end{equation}
provided of course that all sums converge.

Bailey used \eqref{abcd} for proving identities of the 
Rogers--Ramanujan type \cite{Bailey49}.
First he showed that a limiting case of the $q$-Gauss sum 
\cite[Eq. (II.8)]{GR90} provides the following conjugate Bailey pair 
relative to $a$:
\begin{equation}\label{gdr12}
\gamma_L=\frac{a^L q^{L^2}}{(aq)_{\infty}} \quad\text{and}\quad
\delta_L=a^L q^{L^2}
\qquad (|q|<1).
\end{equation}
Substituted into \eqref{abcd} this gives
\begin{equation}\label{BL}
\frac{1}{(aq)_{\infty}}\sum_{L=0}^{\infty}a^L q^{L^2}\alpha_L=
\sum_{L=0}^{\infty} a^L q^{L^2}\beta_L.
\end{equation}
To obtain the Rogers--Ramanujan identities, Bailey now required 
\cite{Bailey47}
\begin{equation}\label{BPR}
\alpha_L=\frac{1-aq^{2L}}{1-a}\:\frac{(-a)^L q^{L(3L-1)/2}(a)_L}
{(q)_L} \quad\text{and}\quad
\beta_L=\frac{1}{(q)_L},
\end{equation}
which follows from a limit of Rogers' $q$-analogue of Dougall's ${_5}F_4$
sum \cite[Eq. (II.21)]{GR90}. Inserting \eqref{BPR} into \eqref{BL} yields
the Rogers--Selberg identity. 
Taking $a=1$ and $a=q$ and performing the sum on the
left using Jacobi's triple product identity \cite[Eq. (II.28)]{GR90}
results in the Rogers--Ramanujan identities \cite{Rogers94}
\begin{equation}\label{RRid}
\sum_{n=0}^{\infty}\frac{q^{n^2}}{(q)_n}=
\frac{1}{(q,q^4;q^5)_{\infty}}\quad\text{ and }\quad
\sum_{n=0}^{\infty}\frac{q^{n(n+1)}}{(q)_n}=
\frac{1}{(q^2,q^3;q^5)_{\infty}},
\end{equation}
where $(a_1,\dots,a_k;q)_r=(a_1;q)_r\cdots(a_k;q)_r$ and $|q|<1$.

About a decade ago Milne and Lilly generalized the notion of a 
Bailey pair to higher dimensions \cite{ML92}. (For a recent review see
\cite{Milne01}.) They defined, for example, an A$_{n-1}$ Bailey pair 
relative to $a$ by
\begin{equation*}
\beta_L=\sum_{\substack{r_i=0 \\ 1\leq i\leq n-1}}^{L_i} M_{L,r} \alpha_r,
\end{equation*}
where $L,r\in\Z_{+}^{n-1}$ and 
\begin{equation*}
M_{L,r}=\prod_{j=1}^{n-1}
\Bigr[(aqx_j/x_{n-1})_{L_j+|r|}\prod_{k=1}^{n-1}
(q^{r_j-r_k+1}x_j/x_k)_{L_j-r_k}\Bigl]^{-1},
\end{equation*}
with $x_1,\dots,x_{n-1}$ indeterminates and $|r|=r_1+\dots+r_{n-1}$.
{}From the point of view of A$_{n-1}$ basic 
hypergeometric series the Milne--Lilly definition of an A$_{n-1}$ 
Bailey pair is not only natural but, as the papers 
\cite{Milne97,ML92,ML95} attest, has also been very fruitful.
However, from the point of view of symmetric functions a rather different 
definition of an A$_{n-1}$ Bailey pair (and conjugate Bailey pair) seems 
natural. To motivate this, let us assume that $a=q^{\ell}$ in \eqref{BP} 
with $\ell$ a nonnegative integer. Then 
\begin{equation}\label{betaL}
\beta_L=\sum_{r=0}^L \frac{\alpha_r}{(aq)_{2L}}\qbin{2L+\ell}{L-r}_q,
\end{equation}
where $\qbins{n}{m}_q$ is a $q$-binomial coefficient, defined as
\begin{equation}\label{qbinomial}
\qbin{n}{m}_q=\begin{cases}
\displaystyle
\frac{(q)_n}{(q)_m(q)_{n-m}} & \text{for } 0\leq m \leq n \\[3mm]
0 & \text{otherwise.}\end{cases}
\end{equation}
For $\lambda,\mu$ partitions and $\eta$ a composition let
$K_{\lambda\eta}$ and $K_{\lambda\mu}(q)$ be a Kostka number
and Kostka polynomial, respectively. (For their definition and for 
notation concerning partitions and compositions, see 
Section~\ref{sec21}.) From \cite[Ch. I, \S 6, Ex. 2]{Macdonald95} 
it follows that for $\lambda=(\lambda_1,\lambda_2)$ and
$\eta=(\eta_1,\eta_2)$ there holds 
\begin{equation}\label{Kchi}
K_{\lambda\eta}=\chi(\lambda_2\leq\eta_1\leq\lambda_1)
\chi(|\lambda|=|\eta|),
\end{equation}
where $\chi(\text{true})=1$ and $\chi(\text{false})=0$.
Also, from \cite[Ch. III, \S 6, Ex. 2]{Macdonald95}
\begin{equation*}
K_{(1^{i-2j} 2^j),(1^i)}(q)=\qbin{i}{j}_q-\qbin{i}{j-1}_q.
\end{equation*}
Combining these formulae leads to the expansion
\begin{align}\label{qbKK}
\qbin{L}{r}_q&=\sum_{k=0}^{\min\{r,L-r\}}
\biggl\{\qbin{L}{k}_q-\qbin{L}{k-1}_q\biggr\} \\
&=\sum_{k=0}^{\lfloor L/2\rfloor} K_{(L-k,k),(r,L-r)}
K_{(1^{L-2k}2^{k}),(1^L)}(q) \notag \\ 
& =\sum_{\mu\vdash L}K_{\mu,(r,L-r)}K_{\mu',(1^L)}(q), \notag
\end{align}
where in the last step we have used that $K_{\mu,(r,L-r)}=0$ if $l(\mu)>2$.
Substituted into \eqref{betaL} this results in
\begin{align}\label{BP2}
\beta_L&=\sum_{r=0}^L \frac{\alpha_r}{(aq)_{2L}}
\sum_{\mu\vdash 2L+\ell}K_{\mu,(L-r,L+r+\ell)}K_{\mu',(1^{2L+\ell})}(q) \\
&=(q)_{\ell}\sum_{r=0}^{(|\eta|-\ell)/2} \frac{\alpha_r}{(q)_{m(\eta)}}
\sum_{\mu\vdash |\eta|}K_{\mu,\frac{|\eta|}{2}(1^2)-k-\ell\bar{\Lambda}_1}
K_{\mu'\eta}(q) \notag
\end{align}
where $\eta=(1^{2L+\ell})$, $(q)_{m(\eta)}=(q)_{2L+\ell}$,
$\bar{\Lambda}_1=(1/2,-1/2)$ and $k=(r,-r)$.
This suggest that a natural multivariable or A$_{n-1}$ generalization
of a Bailey pair is obtained by letting $\eta$ be a partition of
at most $n-1$ different parts and by replacing the composition 
$|\eta|(1^2)/2-k+\ell\bar{\Lambda}_1$ by an appropriate
composition of $n$ parts.
Such a generalization was introduced in \cite{SW99}, and various conjectures
were made concerning A$_{n-1}$ Bailey pairs and conjugate Bailey pairs.
Further progress was made in \cite{ASW99}, where
it was shown that the A$_2$ generalization of \eqref{BP2} leads to 
A$_2$ versions of the Rogers--Ramanujan identities \eqref{RRid} such as
\begin{equation*}
\sum_{n_1,n_2=0}^{\infty}\frac{q^{n_1^2-n_1 n_2+n_2^2}}{(q)_{n_1}}
\qbin{2n_1}{n_2}_q=\frac{1}{(q,q,q^3,q^4,q^6,q^6;q^7)_{\infty}}.
\end{equation*}

In this paper we will study the more general A$_{n-1}$ generalization
of \eqref{BP2} and prove the conjectured conjugate Bailey pairs of 
\cite{SW99}. Unlike in \cite{ASW99,SW98a,SW00b,SW00a,W99},
where conjugate Bailey pairs are proven using summation and transformation
formulas for basic hypergeometric series, our approach here will be
to employ the theory of Kostka polynomials.
Even in the classical or A$_1$ case this is new.
Having proven the conjectures of \cite{SW99}, we will utilize a
conjectural expansion for Kostka numbers to derive much more general
conjugate Bailey pairs in terms of fractional-level string functions
and configuration sums of A$_{n-1}^{(1)}$.
Not only will this give rise to new expressions for admissible
characters of A$_{n-1}^{(1)}$, but also to 
$q$-series identities for unitary as well as nonunitary A-type branching 
functions.

The outline of this paper is as follows.
The next section contains an introduction to symmetric functions,
Kostka polynomials and supernomial coefficients, and concludes with a
conjectured expansion of the Kostka numbers in terms of antisymmetric
supernomials.
In Section~\ref{sec3} we define an A$_{n-1}$ analogue of
Bailey's lemma and prove a generalization of the conjugate
Bailey pair \eqref{gdr12}, settling a conjecture of \cite{SW99}.
Section~\ref{sec4} provides a bosonic reformulation of the
results of Section~\ref{sec3}, which is exploited in the subsequent
section to yield very general conjugate Bailey pairs in terms of 
A$_{n-1}^{(1)}$ string functions and configuration sums.
Finally, in Section~\ref{sec6}, we use these results to find many
new A-type $q$-series identities. 

\section{Symmetric functions, Kostka polynomials 
and supernomial coefficients}\label{sec2}

\subsection{Introduction}\label{sec21}

This section contains a brief introduction to symmetric functions
and Kostka polynomials. For further details, see
\cite{Butler94,Fulton97,Macdonald95,Manivel01}.

A composition $\lambda=(\lambda_1,\lambda_2,\dots)$ is a sequence 
(finite or infinite) of nonnegative integers with finitely many 
$\lambda_i\geq 0$ unequal to zero. 
Two compositions that only differ in their tails of zeros are identified,
and, for example, we will not distinguish between $(1,0,2,3,0,0)$ and 
$(1,0,2,3)$.  Even when writing $\lambda=(\lambda_1,\dots,\lambda_n)
\in\Z_{+}^n$ for a composition $\lambda$ we do give meaning to 
$\lambda_{n+k}$ for $k\geq 1$ (namely the obvious $\lambda_{n+k}=0$).
The weight of a composition $\lambda$, denoted $|\lambda|$, is
the sum of its (nonzero) components. The largest and smallest components of a 
composition $\lambda$ will be denoted by $\max(\lambda)$ and $\min(\lambda)$.
Also for noninteger sequences $\mu$ will we write $|\mu|$ for the sum
of its components.

A partition $\lambda$ is a composition whose components are weakly decreasing.
The nonzero components of a partition $\lambda$ are called its parts
and the number of parts, denoted $l(\lambda)$, its length.
We say that $\lambda$ is a partition of $k$, denoted $\lambda\vdash k$ 
if $|\lambda|=k$. 

The Ferrers graph of a partition $\lambda$ is obtained by drawing
$l(\lambda)$ left-aligned rows of dots, with the $i$th row containing
$\lambda_i$ dots. A partition is called rectangular 
if all its parts are equal, i.e., if its Ferrers graph has rectangular shape.
The conjugate of a partition $\lambda$, denoted $\lambda'$,
is obtained by transposing the Ferrers graph of $\lambda$.
For $\lambda$ a partition,
\begin{equation*}
n(\lambda)=\sum_{i=1}^{l(\lambda)}(i-1)\lambda_i
=\sum_{i=1}^{\lambda_1}\binom{\lambda'_i}{2}.
\end{equation*}
The unique partition $(0,0,\dots)$ of $0$
is denoted by $\emptyset$, the set of all partitions of $k$ by $\Pa_k$ and
the set of all partitions by $\Pa$.
We define the usual dominance (partial) order on $\Pa_k$
by writing $\lambda\geq\mu$ for $\lambda,\mu\in\Pa_k$ iff
$\lambda_1+\cdots+\lambda_i\geq \mu_1+\cdots+\mu_i$ for all $i\geq 1$.
Note that $\lambda\geq\mu\Leftrightarrow\mu'\geq\lambda'$
and that $\lambda\geq\mu$ implies $\lambda_1\geq\mu_1$ as well as
$l(\lambda)\leq l(\mu)$.

Frequently $\lambda\in\Pa$ will be written as
$\lambda=(1^{m_1}2^{m_2}\dots)$ where $m_j=\lambda_j'-\lambda'_{j+1}$ 
is the multiplicity of the part $j$. The associated composition 
$(m_1,m_2,\dots,m_{\lambda_1})$ is denoted by $m(\lambda)$.
Using the multiplicities, the union of the partitions 
$\lambda=(1^{m_1}2^{m_2}\dots)$ and 
$\mu=(1^{n_1}2^{n_2}\dots)$ is given by
$\lambda\cup\mu=(1^{m_1+n_1}2^{m_2+n_2}\dots)$.
For $\lambda$ a partition such that $\lambda_1\leq a$
\begin{equation}\label{nab}
n((a^b)\cup \lambda)=a\binom{b}{2}+b|\lambda|+n(\lambda).
\end{equation}
One can also define the sum of partitions
but, more generally, we will simply assume the standard addition 
and subtraction for sequences in $\mathbb{R}^n$.
For two partitions $\lambda$ and $\mu$ we write $\mu\subseteq\lambda$ if
$\mu_i\leq\lambda_i$ for all $i$, and $\nu=\lambda\cap\mu$ if
$\nu_i=\min\{\lambda_i,\mu_i\}$ for all $i$.

Let $S_n$ be the symmetric group, i.e., the group with elements given by
the permutations of $(1,\dots,n)$, and with multiplication $\circ$
given by the usual composition of permutations.
Given a permutation $\sigma\in S_n$ let $\ell(\sigma)$ be its length 
($=$ the minimal number of ``adjacent transpositions'' required to 
obtain $\sigma$ from $(1,2,\dots,n)$) and $\epsilon(\sigma)
=(-1)^{\ell(\sigma)}$ its signature.
Clearly $\epsilon(\sigma\circ\tau)=\epsilon(\sigma)\epsilon(\tau)$ and (thus)
$\epsilon(\sigma)=\epsilon(\sigma^{-1})$.
The unique permutation of maximal length ($=n(n-1)/2$) will be denoted by
$\pi$. (Incidentally, $\pi$ is the partition $(n,\dots,2,1)$ which also
has length $l(\pi)=n$ in the sense of partitions.)
Let $k=(k_1,\dots,k_n)$. Then $\sigma\in S_n$ acts on $k$ by permuting 
its components; $\sigma(k)=(k_{\sigma_1},\dots,k_{\sigma_n})$.
For $\lambda$ a partition such that $l(\lambda)\leq n$, 
$S_n^{\lambda}$ denotes the subgroup of $S_n$
consisting of permutations that leave $\lambda$ invariant, i.e., 
$\sigma(\lambda)=\lambda$ for $\sigma\in S_n^{\lambda}$.

Let $x^{\lambda}=x_1^{\lambda_1}\dots x_n^{\lambda_n}$. 
For $\lambda$ a partition such that $l(\lambda)\leq n$ the monomial 
symmetric function in $n$ variables is defined as
\begin{equation*}
m_{\lambda}(x)=\sum_{\sigma\in S_n/S_n^{\lambda}} x^{\sigma(\lambda)}.
\end{equation*}
If $l(\lambda)>n$, $m_{\lambda}(x)=0$.
To define the complete symmetric function $h_{\lambda}$, set $h_r(x)=0$ 
for $r<0$, 
\begin{equation*}
h_r(x)=\sum_{\mu\vdash r} m_{\mu}(x)=
\sum_{1\leq i_1\leq\dots\leq i_r\leq n} x_{i_1}\cdots x_{i_r}
\end{equation*}
for $r\geq 0$, and $h_{\lambda}(x)=h_{\lambda_1}(x)h_{\lambda_2}(x)\cdots$
for $\lambda$ a composition.
In much the same way the elementary symmetric function $e_{\lambda}$
is defined by $e_0(x)=1$,
\begin{equation*}
e_r(x)=m_{(1^r)}(x)=\sum_{1\leq i_1<\dots<i_r\leq n} x_{i_1}\cdots x_{i_r}
\end{equation*}
for $r\geq 1$, and $e_{\lambda}(x)=e_{\lambda_1}(x)e_{\lambda_2}(x)\cdots$ 
for $\lambda$ a composition. Since $e_r(x)=0$ for $r>n$ we have 
$e_{\lambda}(x)=0$ if $\max(\lambda)>n$ ($\lambda_1>n$ for $\lambda\in\Pa$).
Clearly, $e_{\lambda}(x)=e_{\sigma(\lambda)}(x)$ and 
$h_{\lambda}(x)=h_{\sigma(\lambda)}(x)$.

For $\lambda\in\Pa$ such that $l(\lambda)\leq n$ let
\begin{equation}\label{a}
a_{\lambda}(x)=\sum_{\sigma\in S_n}\epsilon(\sigma) x^{\sigma(\lambda)}.
\end{equation}
Then the Schur function $s_{\lambda}(x)$ is defined as
\begin{equation*}
s_{\lambda}(x)=\frac{a_{\lambda+\delta}(x)}{a_{\delta}(x)},
\end{equation*}
where $\delta=(n-1,n-2,\dots,1,0)$.
The denominator on the right is known as the Vandermonde determinant;
\begin{equation}\label{Vdet}
a_{\delta}(x)=\prod_{1\leq i<j\leq n}(x_i-x_j).
\end{equation}

Finally we introduce the Hall-Littlewood symmetric function as the following
$q$-analogue of the monomial symmetric function:
\begin{equation*}
P_{\lambda}(x;q)=\sum_{\sigma\in S_n/S_n^{\lambda}} x^{\sigma(\lambda)}
\prod_{\lambda_i>\lambda_j}
\frac{x_{\sigma_i}-qx_{\sigma_j}}{x_{\sigma_i}-x_{\sigma_j}}
\end{equation*}
for $\lambda\in\Pa$ such that $l(\lambda)\leq n$, and
$P_{\lambda}(x;q)=0$ for $\lambda\in\Pa$ such that $l(\lambda)>n$.
Besides the obvious $P_{\lambda}(x;1)=m_{\lambda}(x)$ one also has
$P_{\lambda}(x;0)=s_{\lambda}(x)$.

For $\mu$ a composition and $\lambda\vdash|\mu|$,
the Kostka number $K_{\lambda\mu}$ is defined by
\begin{equation}\label{hs}
h_{\mu}(x)=\sum_{\lambda\vdash|\mu|}K_{\lambda\mu}s_{\lambda}(x).
\end{equation}
More generally, $K_{\lambda\mu}=0$ if $\lambda$ is not a partition,
or $\mu$ not a composition, or $|\mu|\neq |\lambda|$. From 
$h_{\mu}(x)=h_{\sigma(\mu)}(x)$ it follows that
\begin{equation}\label{Ksymm}
K_{\lambda\mu}=K_{\lambda,\sigma(\mu)}.
\end{equation}
Further simple properties of the Kostka numbers are
\begin{equation}\label{Kzero}
K_{\lambda\mu}\neq 0 \text{ iff $\lambda\geq\mu$}
\end{equation}
for $\lambda\in \Pa$ and $\mu\vdash|\lambda|$, and 
\begin{equation}\label{square}
K_{\lambda\mu}=K_{(a^n)+\lambda,(a^n)+\mu}
\end{equation}
for $l(\lambda)\leq n$, $\mu\in\Z_{+}^n$ and $a\geq -\lambda_n$.

Occurrences of the Kostka numbers similar to \eqref{hs} are given by
\begin{equation}\label{es}
e_{\mu}(x)=\sum_{\lambda\vdash|\mu|}K_{\lambda'\mu}s_{\lambda}(x)
\end{equation}
for $\mu$ a composition, and 
\begin{equation}\label{sm}
s_{\lambda}(x)=\sum_{\mu\vdash|\lambda|}K_{\lambda\mu}m_{\mu}(x)
\end{equation}
for $\lambda\in\Pa$. Note that both sides of \eqref{es} vanish
if $\max(\mu)>n$, and both sides of \eqref{sm} vanish if $l(\lambda)>n$.

A useful formula involving the Kostka numbers follows from
the Jacobi-Trudi identity
\begin{equation*}
s_{\mu}(x)=\det(h_{\mu_i+j-i}(x))_{1\leq i,j\leq n}
=\sum_{\sigma\in S_n}\epsilon(\sigma)h_{\sigma(\mu+\delta)-\delta}.
\end{equation*}
Substituting \eqref{hs} and interchanging sums gives
\begin{equation*}
s_{\mu}(x)=\sum_{\lambda\vdash|\mu|}s_{\lambda}(x)
\sum_{\sigma\in S_n}\epsilon(\sigma)K_{\lambda,\sigma(\mu+\delta)-\delta},
\end{equation*}
which implies that
\begin{equation}\label{Kdelta}
\sum_{\sigma\in S_n}\epsilon(\sigma)
K_{\lambda,\sigma(\mu+\delta)-\delta}=\delta_{\lambda,\mu}
\end{equation}
for $\lambda,\mu\in\Pa$ such that $l(\mu)\leq n$ and 
$\delta_{a,b}=\chi(a=b)$. Note that \eqref{Kdelta} also holds
for $\mu=(\mu_1,\mu_2,\dots,\mu_n)$ such that $\mu_1
\geq\mu_2\geq \dots\geq\mu_n$ and $\mu_n<0$, since
$\sigma(\mu+\delta)-\delta$ is not a composition for such a $\mu$.

For $\lambda,\mu\in\Pa$ the Kostka polynomial (or Kostka--Foulkes polynomial) 
$K_{\lambda\mu}(q)$ is defined by 
\begin{equation}\label{sP}
s_{\lambda}(x)=\sum_{\mu\vdash|\lambda|}K_{\lambda\mu}(q)P_{\mu}(x;q)
\end{equation}
and $K_{\lambda\mu}(q)=0$ if $|\lambda|\neq|\mu|$.
We extend this to compositions $\mu$ through $K_{\lambda\mu}(q)=
K_{\lambda,\sigma(\mu)}(q)$ and to more general sequences $\lambda$ and $\mu$
by $K_{\lambda\mu}=0$ if $\lambda$ is not a partition,
or $\mu$ is not a composition.
Since $P_{\lambda}(x;1)=m_{\lambda}(x)$ it follows that
$K_{\lambda\mu}(1)=K_{\lambda\mu}$.

Let $\lambda,\mu,\in\Pa$ such that $r$ is rectangular and
$\lambda,\mu\subseteq r$, and denote by $\tilde{\lambda}_r$ 
($\tilde{\mu}_r$) the ``complement of $\lambda$ ($\mu$) with
respect to $r$''. That is, for $r=(a^b)$,
\begin{equation*}
\tilde{\lambda}_r=(\underbrace{a,\dots,a}_{b-l(\lambda)},
a-\lambda_{l(\lambda)},\dots,a-\lambda_1).
\end{equation*}
Then the following duality holds \cite{LS81}:
\begin{equation}\label{Kcomp}
K_{\lambda\mu}(q)=K_{\tilde{\lambda}_r,\tilde{\mu}_r}(q).
\end{equation}
For example, by taking $r=(4^3)$, 
\begin{equation*}
K_{(3,1),(2,1,1)}(q)=K_{(4,3,1),(3,3,2)}(q)=q+q^2.
\end{equation*}
Note that \eqref{Kcomp} implies a $q$-analogue of \eqref{square}
as well as
\begin{equation}\label{sq2}
K_{\lambda\mu}(q)=K_{(n^a)\cup\lambda,(n^a)\cup\mu}(q)
\end{equation}
for $\lambda,\mu\in\Pa$ such that $\lambda_1,\mu_1\leq n$ and $a\geq 0$.

\subsection{A$_{n-1}$ supernomial coefficients}
We now come to perhaps the two most important definitions,
that of the completely symmetric and
completely antisymmetric A$_{n-1}$ supernomials \cite{HKKOTY99,SW98b,SW99}.
Although all that is needed for their definition is at our disposal,
we shall somewhat digress and motivate the supernomials from
the point of view of generalized Kostka polynomials.
In \cite{KS02,SW99,Shimozono01a,Shimozono01b,Shimozono02,SW00} 
generalizations of the Kostka polynomials $K_{\lambda\mu}(q)$ were
introduced wherein the composition $\mu$ in the ordinary Kostka
polynomials is replaced by a (finite) sequence $R=(R_1,R_2,\dots)$ of 
rectangular partitions $R_i$, such that 
\begin{equation}\label{KtoK}
K_{\lambda R}(q)=K_{\lambda\mu}(q) \quad \text{if $R_i=(\mu_i)$ for all $i$}. 
\end{equation}
(For $R=(R_1)$, $K_{\lambda R}(q)=\delta_{\lambda,R_1}$ 
and \textit{not} $K_{\lambda R}(q)=K_{\lambda R_1}(q)$.)
Whereas $K_{\lambda\mu}$ is the cardinality of the set of (semi-standard) 
Young tableaux of shape $\lambda$ and content (or weight) $\mu$ 
\cite{Macdonald95}, 
$K_{\lambda R}(1)$ is the cardinality of the set of Littlewood--Richardson 
tableaux of shape $\lambda$ and content $R$ \cite{LR34}. 
This implies that $K_{\lambda R}(q)=0$ if $|\lambda|\neq\sum_i|R_i|=:|R|$.
Equivalently, the generalized Kostka polynomials 
may be viewed as $q$-analogues of generalized
Littlewood--Richardson coefficients. Recalling the definition of the
Littlewood--Richard\-son coefficients $c_{\lambda\mu}^{\nu}$ as
\begin{equation}\label{LRc}
s_{\lambda}(x)s_{\mu}(x)=\sum_{\nu\in\Pa}c_{\lambda\mu}^{\nu}s_{\nu}(x)
\end{equation}
for $\lambda,\mu\in\Pa$ (hence $c_{\lambda\mu}=0$ if 
$|\nu|\neq |\lambda|+|\mu|$), there holds
\begin{equation}\label{KR1}
K_{\lambda R}(1)=\sum_{\nu^{(1)},\dots,\nu^{(k-2)}\in\Pa}
c_{\nu^{(0)}R_2}^{\nu^{(1)}}
c_{\nu^{(1)}R_3}^{\nu^{(2)}}\cdots
c_{\nu^{(k-2)}R_k}^{\nu^{(k-1)}},
\end{equation}
with $R=(R_1,\dots,R_k)$, $\nu^{(0)}=R_1$ and $\nu^{(k-1)}=\lambda$.
For $k=2$ this yields $K_{\lambda,(R_1,R_2)}(1) =c_{R_1 R_2}^{\lambda}$.

Now let $\lambda\in\Z_{+}^n$. Then the A$_{n-1}$ supernomials are defined in 
terms of the generalized Kostka polynomials as
\begin{equation}\label{SR}
S_{\lambda R}(q)=\sum_{\nu\vdash|\lambda|}K_{\nu\lambda}K_{\nu' R^{\ast}}(q),
\end{equation}
where $R^{\ast}=(R_1',R_2',\dots)$. From the duality 
\cite[Prop. 12]{KS02}, \cite[Thm. 7.1]{SW99}
\begin{equation}\label{Kdual}
K_{\lambda R}(q)=q^{\sum_{i<j}|R_i\cap R_j|}K_{\lambda' R^{\ast}}(1/q)
\end{equation}
it follows that $S_{\lambda R}(1)=\sum_{\nu\vdash|\lambda|} 
K_{\nu\lambda} K_{\nu R}(1)$. Multiplying this by $m_{\lambda}(x)$,
summing over $\lambda$ and using \eqref{sm}, \eqref{LRc} and \eqref{KR1},
one finds
\begin{equation}\label{sSm}
s_{R_1}(x) s_{R_2}(x)\cdots =\sum_{\lambda\vdash |R|}
S_{\lambda R}(1) m_{\lambda}(x)=
\sum_{\substack{\lambda\in\Z_{+}^n \\|\lambda|=|R|}}
S_{\lambda R}(1) x^{\lambda}.
\end{equation}
The $S_{\lambda R}(1)$ may thus be viewed as a generalized multinomial
coefficient.

In the following we need two special cases of the A$_{n-1}$ supernomials.
One, the completely antisymmetric A$_{n-1}$ supernomial
(antisymmetric supernomial for short) 
$\A_{\lambda\mu}(q)$, arises when $\lambda\in\Z_{+}^n$ and
$R_i=(1^{\mu_i})$ for all $i$.
Since $R_i^{\ast}=R_i'=(\mu_i)$, \eqref{KtoK} and \eqref{SR} imply
\begin{equation}\label{AKK}
\A_{\lambda\mu}(q):=S_{\lambda ((1^{\mu_1}),(1^{\mu_2}),\dots)}(q)
=\sum_{\nu\vdash|\lambda|} K_{\nu\lambda} K_{\nu'\mu}(q).
\end{equation}
By \eqref{es} and \eqref{sP} it thus follows that 
$e_{\mu}(x)=\sum_{\lambda\vdash |\mu|}\A_{\lambda\mu}(q) 
P_{\lambda}(x;q)$.
Since $P_{\lambda}(x;1)=m_{\lambda}(x)$ and
$e_{\mu}(x)=s_{(1^{\mu_1})}(x)s_{(1^{\mu_2})}(x)\cdots$,
this is in accordance with \eqref{sSm}.

In similar fashion the completely symmetric A$_{n-1}$ supernomial
$\Su_{\lambda\mu}(q)$ results when $\lambda\in\Z_{+}^n$,
$R_i=(\mu_i)$ for all $i$, and $q\to 1/q$;
\begin{equation}\label{SKK}
\Su_{\lambda\mu}(q):=
q^{\sum_{i<j}\min\{\mu_i,\mu_j\}}
S_{\lambda ((\mu_1),(\mu_2),\dots)}(1/q)
=\sum_{\nu\vdash|\lambda|} K_{\nu\lambda} K_{\nu\mu}(q),
\end{equation}
where the expression on the right follows from \eqref{KtoK}, \eqref{SR}
and \eqref{Kdual}. Observe that when $\mu$ is a partition
$\sum_{i<j}\min\{\mu_i,\mu_j\}=n(\mu)$.
By \eqref{hs} and \eqref{sP} it follows that
$h_{\mu}(x)=\sum_{\lambda\vdash |\mu|}\Su_{\lambda\mu}(q)
P_{\lambda}(x;q)$.
Since $P_{\lambda}(x;1)=m_{\lambda}(x)$ and
$h_{\mu}(x)=s_{(\mu_1)}(x)s_{(\mu_2)}(x)\cdots$,
this is again in agreement with \eqref{sSm}.

The antisymmetric supernomial of \eqref{AKK} will
play a particularly important role throughout the paper
and it will be convenient to introduce the further notation
\begin{equation}\label{SKK2}
\qbin{L}{\lambda}=\A_{\lambda\mu}(q)=
\sum_{\nu\vdash |\lambda|}K_{\nu\lambda}K_{\nu'\mu}(q),
\end{equation}
for $\lambda\in\Z_{+}^n$, $\mu\in\Pa_{\lambda}$ such that $\mu_1\leq n-1$
and $L=m(\mu)\in\Z_{+}^{n-1}$.
When $\lambda\in\Z^n$ with one or more components being negative we set 
$\qbin{L}{\lambda}=0$.
It is in fact not a restriction to assume that $\mu$ is a partition
with largest part at most $n-1$.
To see this, assume that $\lambda,\mu\in\Pa$ (with $l(\lambda)\leq n$),
which is harmless since
$\A_{\lambda\mu}(q)=\A_{\sigma(\lambda),\tau(\mu)}(q)$.
Then the summand in \eqref{AKK} is nonzero iff $\nu\geq \lambda$ and 
$\nu'\geq\mu$ implying that $\A_{\lambda\mu}(q)$ is nonzero iff 
$\mu'\geq\lambda$. One may therefore certainly assume that
$\mu_1\leq l(\lambda)\leq n$. But from \eqref{square} and \eqref{sq2}
it is easily seen that
$\A_{\lambda+(1^n),(n)\cup \mu}(q)=\A_{\lambda\mu}(q)$,
so that we may actually assume 
$\mu_1\leq n-1$ in $\A_{\lambda\mu}(q)$,
leading naturally to \eqref{SKK2}.

In the case of A$_1$, $\lambda=(\lambda_1,\lambda_2)$, $\mu=(1^{|\lambda|})$
and $L=m(\mu)=(|\lambda|)$, and from the extremes of \eqref{qbKK} it 
follows that
\begin{equation}\label{A1sup}
\qbin{(\lambda_1+\lambda_2)}{(\lambda_1,\lambda_2)}=
\qbin{\lambda_1+\lambda_2}{\lambda_1}_q,
\end{equation}
with on the right the classical $q$-binomial coefficient \eqref{qbinomial}.
Two important symmetries of the antisymmetric supernomials 
needed subsequently are
\begin{equation}\label{Ssymm}
\qbin{L}{\lambda}=\qbin{L}{\sigma(\lambda)}
\end{equation}
and 
\begin{equation}\label{Srev}
\qbin{(L_1,\dots,L_{n-1})}{\lambda}=
\qbin{(L_{n-1},\dots,L_1)}{(|L|^n)-\lambda}.
\end{equation}
We will also use that the Kostka polynomials may be expressed  
in terms of antisymmetric supernomials as
\begin{equation}\label{Ksup}
K_{\lambda\mu}(q)=\sum_{\sigma\in S_n}
\epsilon(\sigma)\qbin{m(\mu)}{\sigma(\lambda'+\delta)-\delta}.
\end{equation}
Since \eqref{Ssymm} follows from \eqref{Ksymm} and \eqref{SKK2}, 
and \eqref{Ksup} from \eqref{Kdelta} and \eqref{SKK2},
we only need to prove \eqref{Srev}.

\begin{proof}[Proof of \eqref{Srev}]
By \eqref{Kcomp}, $\tilde{\lambda'}_{r'}=(\tilde{\lambda}_r)'$
and $|\nu|+|\tilde{\nu}_r|=|r|$, 
\begin{align*}
\A_{\lambda\mu'}(q)&=\sum_{\nu\vdash|\lambda|}K_{\nu\lambda}
K_{\nu'\mu'}(q) 
=\sum_{\nu\vdash|\lambda|}K_{\tilde{\nu}_r,\tilde{\lambda}_r}
K_{\tilde{\nu'}_{r'},\tilde{\mu'}_{r'}}(q) \\
&=\sum_{\nu\vdash|\lambda|}K_{\tilde{\nu}_r,\tilde{\lambda}_r}
K_{(\tilde{\nu}_r)',(\tilde{\mu}_r)'}(q) 
=\sum_{\tilde{\nu}_r\vdash|\tilde{\lambda}_r|}
K_{\tilde{\nu}_r,\tilde{\lambda}_r}K_{(\tilde{\nu}_r)',(\tilde{\mu}_r)'}(q)\\
&=\A_{\tilde{\lambda}_r,(\tilde{\mu}_r)'}(q)
\end{align*}
for $r$ a rectangular partition such that $\lambda,\mu\subseteq r$.
Note that the summand in the second expression of the
first line is nonzero for $\lambda\leq\nu\leq\mu$ only.
This implies that if $\mu,\lambda\subseteq r$ then also $\nu\subseteq r$,
and thus permits the use of \eqref{Kcomp} resulting in the
third expression of the top-line.

Now specialize $\lambda=(\lambda_1,\dots,\lambda_n)\in\Pa$, 
$\mu'=(1^{L_1},\dots,(n-1)^{L_{n-1}})$ and $r=(\mu_1^n)=(|L|^n)$,
where $L=m(\mu')$.
Then $\tilde{\lambda}_r=(|L|^n)-\pi(\lambda)$ and
$(\tilde{\mu}_r)_i=L_1+\cdots+L_{n-i}$ so that
$(\tilde{\mu}_r)'=(1^{L_{n-1}},\dots,(n-1)^{L_1})$.
Recalling \eqref{SKK2}, \eqref{Srev} now follows
for $\lambda\in\Pa$ with the bottom-entry of the right side being
substituted by $(|L|^n)-\pi(\lambda)$. Thanks to \eqref{Ssymm}
this may be replaced by $(|L|^n)-\lambda$ with $\lambda$ a composition. 
\end{proof}

\subsection{A conjecture}
To conclude this section we present a conjectured expansion of the
Kostka numbers in terms of antisymmetric supernomials.
Admitting this conjecture leads to an orthogonality relation between
Kostka polynomials and antisymmetric supernomials (Corollary~\ref{corinv} 
below), which will be crucial in Section~\ref{sec5}.

Before we can state the conjecture some more notation is needed.
Denote the set $\{1,\dots,k\}$ by $[k]$.
Let $\varepsilon_i$ for $i\in[n]$ be the canonical basis vectors 
in $\mathbb{R}^n$ with standard inner product 
$\ip{\varepsilon_i}{\varepsilon_j}=\delta_{i,j}$. Occasionally we will
abbreviate $\ip{v}{v}=||v||^2$.
Expressed in terms of the $\varepsilon_i$, the simple roots and 
fundamental weights $\alpha_i$ and $\bar{\Lambda}_i$ 
$(i\in[n-1])$ of A$_{n-1}$ are given by
$\alpha_i=\varepsilon_i-\varepsilon_{i+1}$ and 
$\bar{\Lambda}_i=\varepsilon_1+\cdots+\varepsilon_i-\frac{i}{n}(1^n)$,
respectively.
(The notation $\Lambda_i$ will be reserved for the
fundamental weights of A$_{n-1}^{(1)}$.)
The Weyl vector $\rho$ is given by the sum of the fundamental weights;
$\rho=\frac{1}{2}(n-1,n-3,\dots,1-n)=\delta-(n-1)(1^n)/2$.
The A$_{n-1}$ weight and root lattices $P$ and $Q\subset P$ are 
the integral span of the fundamental weights and the integral span of 
the simple roots, respectively.
$P_{+}\subset P$ denotes the set of dominant weights, i.e., those
weights $\lambda\in P$ for which $\ip{\alpha_i}{\lambda}\geq 0$ for all $i$.
Specifically, for $\lambda\in P$ we have
$\lambda\in\Z^n/n$, $|\lambda|=0$ and $\lambda_i-\lambda_{i+1}\in\Z$.
For $\lambda\in P_{+}$ the third condition needs to
be sharpened to $\lambda_i-\lambda_{i+1}\in\Z_{+}$, and for $\lambda\in Q$
the first condition needs to be replaced by $\lambda\in\Z^n$.
The weight and root lattices are invariant under the action of $S_n$,
and for $\lambda\in P$ and $\sigma\in S_n$, 
$\sigma(\lambda)-\lambda\in Q$. This last fact easily follows from
$\sigma^{(j)}(\bar{\Lambda}_i)=\bar{\Lambda}_i-\alpha_i \delta_{i,j}$,
where $\sigma^{(j)}=(1,2,\dots,n)+\alpha_j$ is the $j$th adjacent 
transposition.
If $C$ and $C^{-1}$ are the Cartan and inverse Cartan matrices
of A$_{n-1}$, i.e., $C_{i,j}=2\delta_{i,j}-\delta_{i,j-1}-\delta_{i,j+1}$
and $C^{-1}_{i,j}=\min\{i,j\}-ij/n$, then
$\alpha_i=\sum_{j=1}^{n-1}C_{i,j}\bar{\Lambda}_j$,
$\ip{\bar{\Lambda}_i}{\alpha_j}=\delta_{i,j}$, 
$\ip{\alpha_i}{\alpha_j}=2C_{i,j}$ and 
$\ip{\bar{\Lambda}_i}{\bar{\Lambda}_j}=C^{-1}_{i,j}$.

Finally, introducing the notation 
$(q)_{\lambda}=\prod_{i\geq 1}(q)_{\lambda_i}$ for $\lambda$ a composition,
and using the shorthand $a\equiv b\pod{c}$ for $a\equiv b\pmod{c}$,
our conjecture can be stated as follows.
\begin{conjecture}\label{conj1}
For $\mu\in\Pa$ and $\nu\in\Z^n$ such that
$l(\mu)\leq n$ and $|\mu|=|\nu|$,
\begin{multline*}
K_{\mu\nu}=
\sum_{\substack{\eta\in\Pa\\ \eta_1\leq n-1\\[0.3mm]|\eta|\equiv|\mu|\pmd{n}}}
\sum_{\sigma\in S_n}\sum_{\lambda\in nQ+\sigma(\rho)-\rho}
\epsilon(\sigma)q^{\frac{1}{2n}\ip{\lambda}{\lambda+2\rho}}  \\
\times
\frac{1}{(q)_{m(\eta)}}
\qbin{m(\eta)}{\frac{|\eta|-|\nu|}{n}(1^n)+\nu}
\qbin{m(\eta)}{\frac{|\eta|-|\mu|}{n}(1^n)+\mu-\lambda}.
\end{multline*}
\end{conjecture}
Since the right side satisfies the periodicity 
$f_{\mu\nu}=f_{\mu+(a^n),\nu+(b^n)}$ for $a,b\in\Z$,
the condition $|\mu|=|\nu|$ cannot be dropped.
We also remark that in the notation of \cite{Macdonald95} 
$(q)_{m(\eta)}=b_{\eta}(q)$.

In the following we give proofs of Conjecture~\ref{conj1} for $n=2$
and for $q=0$.
\begin{proof}[Proof for $n=2$]
Set $\eta=(1^r)$ (so that $m(\eta)=|\eta|=r$) and $\lambda=(s,-s)$.
When $\sigma=(1,2)$ ($\sigma=(2,1)$) we need to sum $s$ 
over the even (odd) integers. Using \eqref{A1sup}
the $n=2$ case of the conjecture thus becomes
\begin{equation*}
K_{\mu\nu}=
\sum_{\substack{r=0 \\ r\equiv|\mu|\pmd{2}}}^{\infty}
\sum_{s=-\infty}^{\infty}\frac{(-1)^s q^{\binom{s+1}{2}}}{(q)_r}
\qbin{r}{\frac{1}{2}(r+\nu_{12})}_q
\qbin{r}{\frac{1}{2}(r+\mu_{12})-s}_q,
\end{equation*}
where $\mu_{12}=\mu_1-\mu_2$, $\nu_{12}=\nu_1-\nu_2$
with $\mu=(\mu_1,\mu_2)$ and $\nu=(\nu_1,\nu_2)$
such that $|\mu|=|\nu|$.

Without loss of generality we may assume that $\nu_{12}\geq 0$.
Then, by \eqref{Kchi}, the left side is nothing but 
$\chi(\mu_2\leq\nu_1\leq \mu_1)=\chi(\mu_{12}\geq \nu_{12})$.
This is also true if $\nu_2<0$ as it implies
$K_{\mu\nu}=0$ as well as $\nu_{12}>|\nu|=|\mu|\geq \mu_{12}$.

On the right we now change $s\to (r+\mu_{12})/2-s$ and perform
the sum over $s$ using the $q$-binomial theorem \cite[Thm. 3.3]{Andrews76}
\begin{equation}\label{qbt}
\sum_{k=0}^n(-z)^k q^{\binom{k}{2}}\qbin{n}{k}_q=(z)_n.
\end{equation}
As a result
\begin{multline*}
\chi(\mu_{12}\geq\nu_{12})\\
=\sum_{\substack{r=\nu_{12} \\ r\equiv|\mu|\pmd{2}}}^{\mu_{12}}
(-1)^{(r+\mu_{12})/2} q^{\binom{(r+\mu_{12})/2+1}{2}}
\frac{(q^{-(r+\mu_{12})/2})_r}{(q)_r}
\qbin{r}{\frac{1}{2}(r+\nu_{12})}_q.
\end{multline*}
Replacing $r\to 2r+\nu_{12}$ and making some simplifications yields
\begin{multline*}
\chi(\mu_{12}\geq\nu_{12})=
(-1)^{(\mu_{12}-\nu_{12})/2} q^{\binom{(\nu_{12}-\mu_{12})/2}{2}}
\qbin{\frac{1}{2}(\mu_{12}+\nu_{12})}{\nu_{12}}_q \\
\times \sum_{r=0}^{(\mu_{12}-\nu_{12})/2}
\frac{(q^{(\mu_{12}+\nu_{12})/2+1},q^{-(\mu_{12}-\nu_{12})/2})_r}
{(q,q^{\nu_{12}+1})_r}.
\end{multline*}
By the $q$-Chu--Vandermonde sum \cite[Eq. (II.7)]{GR90}
\begin{equation}\label{qCV}
{_2\phi_1}(a,q^{-n};c;q,cq^n/a):=
\sum_{k=0}^n \frac{(a,q^{-n})_k}{(q,c)_k}\Bigl(\frac{cq^n}{a}\Bigr)^k=
\frac{(c/a)_n}{(c)_n}
\end{equation}
this is readily found to be true.
\end{proof}

\begin{proof}[Proof for $q=0$] 
Replace $\lambda\to\sigma(\lambda+\rho)-\rho=
\sigma(\lambda+\delta)-\delta$ in the sum on the right. This
leaves $\ip{\lambda}{\lambda+2\rho}$ unchanged, yielding
\begin{multline}\label{lambdashift}
K_{\mu\nu}=
\sum_{\substack{\eta\in\Pa\\ \eta_1\leq n-1\\[0.3mm]|\eta|\equiv|\mu|\pmd{n}}}
\sum_{\sigma\in S_n}\sum_{\lambda\in nQ}
\epsilon(\sigma) q^{\frac{1}{2n}\ip{\lambda}{\lambda+2\rho}}  \\
\times
\frac{1}{(q)_{m(\eta)}}
\qbin{m(\eta)}{\frac{|\eta|-|\nu|}{n}(1^n)+\nu}
\qbin{m(\eta)}{\frac{|\eta|-|\mu|}{n}(1^n)+\mu-\sigma(\lambda+\delta)+\delta}.
\end{multline}
Now use that (i)
$q^{\ip{\lambda}{\lambda+2\rho}}=\delta_{\lambda,\emptyset}+O(q)$
for $\lambda\in nQ$, (ii) $\qbin{m(\eta)}{\mu}=K_{\eta'\mu}+O(q)$,
as follows from $K_{\lambda\mu}(0)=\delta_{\lambda,\mu}$ and 
equation \eqref{SKK2}, and (iii) $1/(q)_{m(\eta)}=1+O(q)$.
Hence the constant term of \eqref{lambdashift} can be extracted as
\begin{equation*}
K_{\mu\nu}=
\sum_{\substack{\eta\in\Pa\\ \eta_1\leq n-1\\[0.3mm]|\eta|\equiv|\mu|\pmd{n}}}
\sum_{\sigma\in S_n}\epsilon(\sigma)
K_{\eta',\frac{|\eta|-|\nu|}{n}(1^n)+\nu}
K_{\eta',\frac{|\eta|-|\mu|}{n}(1^n)+\mu-\sigma(\delta)+\delta}.
\end{equation*}
Making the replacement $\sigma\to\sigma^{-1}$ and using the symmetry
\eqref{Ksymm} results in the sum $\sum_{\sigma\in S_n} \epsilon(\sigma)
K_{\eta',\sigma(\lambda+\delta)-\delta}$
with $\lambda=\frac{|\eta|-|\mu|}{n}(1^n)+\mu$.
By \eqref{Kdelta} (plus subsequent comment) and the substitution
$\eta\to\eta'$ this gives
\begin{equation*}
K_{\mu\nu}=
\sum_{\substack{\eta\in\Pa\\l(\eta)\leq n-1\\[0.3mm]|\eta|\equiv|\mu|\pmd{n}}}
K_{\eta,\frac{|\eta|-|\nu|}{n}(1^n)+\nu}
\delta_{\eta,\frac{|\eta|-|\mu|}{n}(1^n)+\mu}.
\end{equation*}
For $\eta,\mu\in\Pa$ such that $l(\eta),l(\mu)\leq n$,
the solution to the linear equation 
$\eta=\frac{|\eta|-|\mu|}{n}(1^n)+\mu$ is
$\eta=\mu+(a^n)$, with $a$ an integer such that 
$-\mu_n\leq a\leq \eta_n$. Since $\eta_n=0$, however,
$a$ is fixed to $a=-\mu_n$. Hence
$\delta_{\eta,\frac{|\eta|-|\mu|}{n}(1^n)+\mu}=
\delta_{\eta,\mu-\mu_n(1^n)}$, and after performing the sum over $\eta$ 
and using that $|\mu|=|\nu|$
the right-hand side becomes $K_{\mu-\mu_n(1^n),\nu-\mu_n(1^n)}$.
Thanks to \eqref{square} this is $K_{\mu\nu}$. 
\end{proof}

The above proof suggests that what is needed to prove
Conjecture~\ref{conj1} is a generalization of \eqref{Kdelta} for
the sum
\begin{equation}\label{l1gen}
\sum_{\sigma\in S_n}
\epsilon(\sigma) 
K_{\lambda,\sigma(\mu+\delta)-\delta-\alpha}
\end{equation}
with $\alpha\in Q$.
Indeed, substituting \eqref{SKK2} into \eqref{lambdashift} we obtain
\begin{multline*}
K_{\mu\nu}=
\sum_{\substack{\eta\in\Pa\\ \eta_1\leq n-1\\[0.3mm]|\eta|\equiv|\mu|\pmd{n}}}
\sum_{\omega\vdash|\eta|}
\frac{K_{\omega'\eta}(q)}{(q)_{m(\eta)}}
\qbin{m(\eta)}{\frac{|\eta|-|\mu|}{n}(1^n)+\nu} \\
\times \sum_{\lambda\in nQ}q^{\frac{1}{2n}\ip{\lambda}{\lambda+2\rho}}
\sum_{\sigma\in S_n}
\epsilon(\sigma) 
K_{\omega,\frac{|\eta|-|\mu|}{n}(1^n)+\mu-\sigma(\lambda+\delta)+\delta}.
\end{multline*}
Changing $\sigma\to\sigma^{-1}$ and using \eqref{Ksymm} the sum over $\sigma$ 
takes exactly the form of \eqref{l1gen} 
(with $\mu\to \frac{|\eta|-|\mu|}{n}(1^n)+\mu$ and $\alpha\to\lambda$).

We also remark that the proof of the constant term cannot (easily) be 
extended to deal with low-order terms in $q$. This because 
$\ip{\lambda}{\lambda+2\rho}=2n$ for 
$\lambda=-n(\bar{\Lambda}_1+\bar{\Lambda}_{n-1})$, 
so that restricting the sum over $\lambda$ to the single term 
$\lambda=\emptyset$ is only correct to zeroth order.

Assuming Conjecture~\ref{conj1} it is not hard to prove the following 
orthogonality relation.
\begin{corollary}\label{corinv}
For $\mu,\nu\in\Pa$ such that $l(\mu),l(\nu)\leq n$ and $|\mu|=|\nu|$,
\begin{multline*}
\sum_{\substack{\eta\in\Pa\\ \eta_1\leq n-1\\[0.3mm]|\eta|\equiv|\mu|\pmd{n}}}
\sum_{\sigma\in S_n}\sum_{\lambda\in nQ+\sigma(\rho)-\rho}
\epsilon(\sigma)q^{\frac{1}{2n}\ip{\lambda}{\lambda+2\rho}}  \\
\times
\frac{K_{\bigl(\frac{|\eta|-|\nu|}{n}(1^n)+\nu\bigr)',\eta}(q)}{(q)_{m(\eta)}}
\qbin{m(\eta)}{\frac{|\eta|-|\mu|}{n}(1^n)+\mu-\lambda}
=\delta_{\mu,\nu}.
\end{multline*}
\end{corollary}
Since $K_{\lambda\mu}(q)=0$ if $\lambda\not\in\Pa$ one may add the 
additional restriction $|\eta|\geq |\nu|-n\nu_n=
\sum_{i=1}^{n-1}(\nu_i-\nu_{i+1})$ to the sum over $\eta$.
We also note that care should be taken when writing
$\bigl(\frac{|\eta|-|\nu|}{n}(1^n)+\nu\bigr)'$ as
$(n^{(|\eta|-|\nu|)/n})\cup\nu'$ since $|\eta|-|\nu|$ need not be nonnegative.

\begin{proof}
Replace $\nu\to \tau(\hat{\nu}+\delta)-\delta$ in Conjecture~\ref{conj1}
with $\hat{\nu}\in\Pa$. Then multiply the result by $\epsilon(\tau)$
and sum $\tau$ over $S_n$. By \eqref{Kdelta} and \eqref{Ksup}  
Corollary~\ref{corinv} (with $\nu$ replaced by $\hat{\nu}$) follows. 
\end{proof}

A proof very similar to the proof of the $n=2$ case of Conjecture~\ref{conj1}
shows that Corollary~\ref{corinv} may be viewed as an A$_{n-1}$ 
generalization of the well-known ${_2\phi_1}(a,q^{-n};aq^{1-n};q,q)=
\delta_{n,0}$, which corresponds to the specialization $c=aq^{1-n}$ in
the $q$-Chu--Vandermonde sum \eqref{qCV}.
Indeed, if we take $n=2$ in Corollary~\ref{corinv} and sum over 
$\sigma$ and $\lambda$ by the $q$-binomial theorem \eqref{qbt}, then 
the remaining sum over $\eta$ takes the form of the above ${_2\phi_1}$ 
sum with $n=\mu_1-\nu_1$ and $a=q^{\mu_1-\nu_2+1}$.

\section{An A$_{n-1}$ Bailey lemma}\label{sec3}

\subsection{Definitions and main result}
Let $\ell\in\Z_{+}$, $k\in Q\cap P_{+}$ (i.e., $k=(k_1,\dots,k_n)\in\Z_{+}^n$
such that $k_1\geq k_2\geq\dots \geq k_n$, and $|k|=0$) and $\eta\in\Pa$ 
such that $\eta_1\leq n-1$ and $|\eta|\equiv \ell\pmod{n}$.
Then $\alpha=(\alpha_k)$ and $\beta=(\beta_{\eta})$ form an A$_{n-1}$ 
Bailey pair relative to $q^{\ell}$ if
\begin{equation}\label{abdef}
\beta_{\eta}=\sum_{\substack{k\in Q\cap P_{+} \\ k_1\leq (|\eta|-\ell)/n}}
\frac{\alpha_k}{(q)_{m(\eta)}}
\qbin{m(\eta)}{\frac{|\eta|}{n}(1^n)-k-\ell\bar{\Lambda}_{n-1}},
\end{equation}
and $\gamma=(\gamma_k)$ and $\delta=(\delta_{\eta})$ form an A$_{n-1}$
conjugate Bailey pair relative to $q^{\ell}$ if
\begin{equation}\label{gddef}
\gamma_k=
\sum_{\substack{\eta\in\Pa\\ \eta_1\leq n-1\\[0.3mm]
|\eta|\geqq \ell+nk_1\pmd{n}}}
\frac{\delta_{\eta}}{(q)_{m(\eta)}}
\qbin{m(\eta)}{\frac{|\eta|}{n}(1^n)-k-\ell\bar{\Lambda}_{n-1}}.
\end{equation}
Here, and elsewhere in the paper, $a\geqq b\pmod{n}$ stands for 
$a\equiv b\pmod{n}$ such that $a\geq b$. 
Note that a necessary condition for the above summands to be non-vanishing 
is $\frac{|\eta|}{n}(1^n)-k-\ell\bar{\Lambda}_{n-1}\in\Z^n_{+}$.
Since $k\in Q\cap P_{+}$ this boils down to the single condition 
$(|\eta|-\ell)/n-k_1\geq 0$, justifying the restrictions in the sums over
$k$ and $\eta$.

Recalling \eqref{SKK2}, the definition \eqref{abdef} is a multivariable 
generalization of \eqref{BP2} up to an irrelevant normalization factor
$(q)_{\ell}$. Of course, to really reduce to \eqref{BP2} for $n=2$
we have to write $\eta=(1^{2L+\ell})$ and then replace 
$\beta_{(1^{2L+\ell})}$ by $\beta_L$.
In exactly the same way \eqref{gddef} is a (normalized) 
multivariable generalization of \eqref{CBP}.

Given an A$_{n-1}$ Bailey pair and conjugate Bailey pair
relative to $q^{\ell}$, it follows that
\begin{equation}\label{abcdgen}
\sum_{k\in Q\cap P_{+}}\alpha_k\gamma_k=
\sum_{\substack{\eta\in\Pa\\ \eta_1\leq n-1\\[0.3mm]|\eta|\geqq\ell\pmd{n}}}
\beta_{\eta}\delta_{\eta}.
\end{equation}
If for $n=2$ we set $\eta=(1^{2L+\ell})$ and write $\beta_L$ and $\delta_L$
instead of $\beta_{\eta}$ and $\delta_{\eta}$
the right side simplifies to $\sum_{L\geq 0}\beta_L\delta_L$ 
in accordance with \eqref{abcd}.

The main result of this section will be an A$_{n-1}$, level-$N$
generalization of Bailey's conjugate pair \eqref{gdr12}, resulting in
an according generalization of Bailey's identity \eqref{BL}. 
First, however, a few more definitions are required.

Let $N$ be a fixed, nonnegative integer, which will be referred to as 
the level. 
The canonical basis vectors in $\mathbb{R}^{n-1}$ and $\mathbb{R}^{N-1}$
will be denoted by $e_a$ and $\bar{e}_j$, respectively, and we set
$\bar{e}_0=\bar{e}_N=\emptyset$.
The $(n-1)$- and $(N-1)$-dimensional identity matrices are 
denoted by $I$ and $\bar{I}$ and 
the A$_{N-1}$ Cartan matrix by $\bar{C}$.
Hence $\bar{C}^{-1}_{j,k}=\min\{j,k\}-jk/N$, which will
also be used when either $j$ or $k$ is $0$ or $N$.
For $u,v\in\Z^k$ with $k$ either $n-1$, $N-1$ or $(n-1)(N-1)$ and
$M$ an $k$-dimensional (square) matrix,
$Mu$ is the vector with components $(Mu)_i=\sum_{j=1}^k M_{i,j}u_j$
and $vMu=\sum_{i,j=1}^k v_i M_{i,j} u_j$.
We already defined $(q)_u=\prod_{i=1}^k (q)_{u_k}$ and extending this,
we also use $\qbins{u+v}{u}_q=\prod_{i=1}^k \qbins{u_i+k_i}{u_i}_q$.
Given integers $m_j^{(a)}$ for $a\in[n-1]$ and $j\in[N-1]$,
the vector $m$ is defined as
\begin{equation}\label{m}
m=\sum_{a=1}^{n-1}\sum_{j=1}^{N-1}(e_a\otimes \bar{e}_j)m_j^{(a)}
 \in\Z^{(n-1)(N-1)}.
\end{equation}

Let $\eta,\nu\in\Pa$ such that $\eta_1\leq n-1$, $\nu_1\leq N$ and 
$|\eta|\equiv|\nu| \pmod{n}$, and set 
$\mu=\sum_{a=1}^{n-1}m_a(\eta)\bar{\Lambda}_a
=\sum_{a=1}^{n-1}(\eta_a'-\eta_{a+1}')\bar{\Lambda}_a\in P_{+}$
($\eta_n'=0$), so that $|\nu|\bar{\Lambda}_1-\mu\in Q$.
Then the polynomial $F_{\eta\nu}(q)$ is defined by
\begin{multline}\label{F}
F_{\eta\nu}(q)=
q^{n(\nu)-\frac{(n-1)|\nu|^2}{2nN}+\frac{\ip{\mu}{\mu}}{2N}} \\
\times \sum_m q^{\frac{1}{2}m (C\otimes\bar{C}^{-1})m
-\sum_{i=1}^{l(\nu)}(e_1\otimes \bar{e}_{N-\nu_i})(I\otimes\bar{C}^{-1})m}
\qbin{m+p}{m}_q.
\end{multline}
Here the vector $p\in\Z^{(n-1)(N-1)}$ is determined by
\begin{equation}\label{p1}
(C\otimes \bar{I})m+
(I\otimes\bar{C})p
=(m(\eta)\otimes\bar{e}_1)+\sum_{i=1}^{l(\nu)}(e_1\otimes\bar{e}_{N-\nu_i})
\end{equation}
and the sum is over $m$ such that
\begin{equation}\label{rest1}
\sum_{a=1}^{n-1}\sum_{j=1}^{N-1}jm_j^{(a)}\alpha_a
\in NQ+\mu-|\nu|\bar{\Lambda}_1.
\end{equation}
This last equation ensures that the components of $p$
are integers, but the reverse of this is not true;
demanding that $p_j^{(a)}\in\Z$ leads to a restriction on $m$ 
that, generally, is weaker than \eqref{rest1}.
We also caution the reader not to confuse $m(\eta)\in\Z_{+}^{n-1}$
with the vector $m\in\Z^{(n-1)(N-1)}$.

For our next definition it will be convenient to view $\mu$ 
rather than $\eta$ as primary variable, and for
$\mu\in P_{+}$ and $\nu\in\Pa$ such that $\nu_1\leq N$ and 
$|\nu|\bar{\Lambda}_1-\mu\in Q$ we define 
\begin{multline}\label{C}
C_{\mu\nu}(q)=\frac{q^{n(\nu)-\frac{(n-1)|\nu|^2}{2nN}
+\frac{\ip{\mu}{\mu}}{2N}}}{(q)_{\infty}^{n-1}} \\
\times
\sum_m\frac{q^{\frac{1}{2}m (C\otimes\bar{C}^{-1})m
-\sum_{i=1}^{l(\nu)}(e_1\otimes \bar{e}_{N-\nu_i})
(I\otimes\bar{C}^{-1})m}}{(q)_m},
\end{multline}
where the sum is again over $m$ such that \eqref{rest1} holds.

{}From \eqref{nab} it readily follows that both $F_{\eta\nu}$
and $C_{\mu\nu}$ satisfy the equation $g_{\xi,(N^n)\cup \nu}(q)=
q^{|\nu|}g_{\xi\nu}(q)$. Without loss of generality we may therefore 
assume that $\nu$ has at most $n-1$ parts equal to $N$.
{}From \eqref{m} it also follows that $F_{\eta\nu}$ and
$C_{\mu\nu}$ trivialize for $N=1$ to
\begin{equation}\label{FN1}
F_{\eta,(1^i)}(q)
=q^{\frac{1}{2}\ip{\mu}{\mu}-\frac{1}{2}||\bar{\Lambda}_i||^2}
\end{equation}
with $\mu=\sum_{a=1}^{n-1}m_a(\eta)\bar{\Lambda}_a$ and
\begin{equation}\label{CN1}
C_{\mu,(1^i)}(q)
=\frac{q^{\frac{1}{2}\ip{\mu}{\mu}-\frac{1}{2}||\bar{\Lambda}_i||^2}}
{(q)_{\infty}^{n-1}},
\end{equation}
where $i\in\{0,\dots,n-1\}$.
For later use it will also be convenient to define
\begin{equation}\label{FCN0}
F_{\eta,\emptyset}(q)=\delta_{\eta,\emptyset}
\quad\text{and}\quad
C_{\mu,\emptyset}(q)=\delta_{\mu,\emptyset} \qquad \text{for $N=0$.}
\end{equation}

In the important special case $\nu=(sN^r)$, the polynomial $F_{\eta\nu}(q)$
can be identified with a level-$N$ restricted version of the generalized 
Kostka polynomials $K_{\lambda R}(q)$ \cite{SS01}.
Denoting these polynomials by $K_{\lambda R}^N(q)$
it follows from \cite[Eq. (6.7)]{SS01} that for
$\nu=(sN^r)$ with $s\in\{0,\dots,N\}$
\begin{equation*}
F_{\eta\nu}(q)=q^{\sum_{a=1}^{n-1}\binom{\eta'_a}{2}
+s+\frac{n-r}{n}(|\eta|-|\nu|)}K^N_{\lambda R}(q^{-1}),
\end{equation*}
where $\lambda=\frac{1}{n}(nN+|\eta|-|\nu|)(1^n)+(s)$ and
$R$ is a sequence of $|m(\eta)|+n-r$ partitions with $m_a(\eta)$ partitions
$(1^a)$ ($a\in[n-1]$) and $n-r$ partitions $(N)$.

Finally we state our main result, which for $n\geq 3$ was 
conjectured in \cite[Eq. (9.9)]{SW99} (see also \cite[Conj. 14]{W99}).
For $n=2$ the theorem below is equivalent to the 
``higher-level Bailey lemma'' of \cite[Cor. 4.1]{SW98a},
and for $N=1$ and $\ell=0$ (but general $n$)
it was proven in \cite[Prop. 12]{W99}.
\begin{theorem}\label{thm1}
Let $\ell\in\Z_{+}$, $k\in Q\cap P_{+}$ and $\eta,\nu\in\Pa$ such that 
$\nu_1\leq N$, $\eta_1\leq n-1$ and $|\eta|\equiv|\nu|\equiv \ell\pmod{n}$.
Then 
\begin{equation*}
\gamma_k=C_{\ell\bar{\Lambda}_1-\pi(k),\nu}(q)\quad\text{and}\quad
\delta_{\eta}=F_{\eta\mu}(q)
\end{equation*}
form an \textup{A}$_{n-1}$ conjugate Bailey pair relative to $q^{\ell}$.
\end{theorem}

Substituting the conjugate Bailey pair of Theorem~\ref{thm1} into
\eqref{abcdgen} we obtain the following generalization of Bailey's 
key-identity \eqref{BL}:
\begin{equation}\label{BLn}
\sum_{k\in Q\cap P_{+}}\alpha_k C_{\ell\bar{\Lambda}_1-\pi(k),\nu}(q)=
\sum_{\substack{\eta\in\Pa\\ \eta_1\leq n-1\\[0.3mm]|\eta|\geqq\ell\pmd{n}}}
\beta_{\eta}F_{\eta\nu}(q).
\end{equation}
Equation \eqref{BL} with $a=q^{\ell}$ follows for $n=2$ and $N=1$. 
According to \eqref{FN1} and \eqref{CN1} one then finds
$C_{(\ell/2+L,\ell/2-L),(1^i)}(q)=
g_{i,\ell}(q) q^{L(L+\ell)}/(q)_{\infty}$ and
$F_{(1^{2L+\ell}),(1^i)}(q)=g_{i,\ell}(q) q^{L(L+\ell)}$,
where $g_{i,\ell}(q)=q^{i(i-2)/4+\ell^2/4}$ is an irrelevant
factor that drops out of \eqref{BLn}.

Applications of \eqref{BLn} will be given in Section~\ref{sec6}, with 
the remainder of this section devoted to a proof of Theorem~\ref{thm1}.
This proof consists of several steps.
First, in Lemma~\ref{lem1} below, we reformulate the 
definition of an A$_{n-1}$ conjugate Bailey pair. 
Then we state Theorem~\ref{thm2}, which claims a deep identity for
Kostka polynomials. A straightforward specialization 
leads to \eqref{corK}, which in combination with Lemma~\ref{lem1} 
proves Theorem~\ref{thm1}. 
To prove Theorem~\ref{thm2} we need the Propositions~\ref{prop1} 
and \ref{prop2} below, which are due to Hatayama \textit{et al.} 
\cite{HKKOTY99} (see also \cite{HKOTY99b}). Because of their importance
to Theorem~\ref{thm2} we have included Sections~\ref{sec33} and \ref{sec34}
containing proofs of both propositions.

\subsection{Proof of Theorem~\ref{thm1}}
\begin{lemma}\label{lem1}
The definition \eqref{gddef} of an A$_{n-1}$ conjugate Bailey pair
\textup{(}relative to $q^{\ell}$\textup{)} can be rewritten as
\begin{equation*}
\gamma_k=\sum_{\substack{\lambda\in\Pa\\ l(\lambda)\leq n-1 \\[0.3mm]
|\lambda|\geqq\ell+nk_1\pmd{n}}}
K_{\lambda,\frac{|\lambda|}{n}(1^n)-k-\ell\bar{\Lambda}_{n-1}}
\sum_{\substack{\eta\in\Pa\\ \eta_1\leq n-1 \\[0.3mm]
|\eta|\geqq|\lambda|\pmd{n}}}
\frac{\delta_{\eta}K_{(n^{(|\eta|-|\lambda|)/n})\cup\lambda',\eta}(q)}
{(q)_{m(\eta)}}.
\end{equation*}
\end{lemma}

\begin{proof}
Using \eqref{SKK2} to eliminate the antisymmetric supernomial from 
\eqref{gddef} yields
\begin{equation*}
\gamma_k=\sum_{\substack{\eta\in\Pa\\ \eta_1\leq n-1\\[0.3mm]
|\eta|\geqq\ell+nk_1\pmd{n}}}
\frac{\delta_{\eta}}{(q)_{m(\eta)}}\sum_{\nu\vdash |\eta|}
K_{\nu,\frac{|\eta|}{n}(1^n)-k-\ell\bar{\Lambda}_{n-1}} K_{\nu'\eta}(q).
\end{equation*}
By \eqref{Ksymm}, \eqref{Kzero} and 
$\min(|\eta|(1^n)/n-k-\ell\bar{\Lambda}_{n-1})=
(|\eta|-\ell)/n-k_1$ (since $k\in Q\cap P_{+}$)
we may add the conditions
$l(\nu)\leq n$ and $(|\eta|-\ell)/n-k_1-\nu_n\in\Z_{+}$ to the sum 
over $\nu$.
Changing the order of summation then gives
\begin{equation*}
\gamma_k=\sum_{\substack{\nu\in\Pa\\l(\nu)\leq n\\[0.3mm]
|\nu|-n\nu_n\geqq\ell+nk_1\pmd{n}}} 
K_{\nu,\frac{|\nu|}{n}(1^n)-k-\ell\bar{\Lambda}_{n-1}} 
\sum_{\substack{\eta\vdash |\nu|\\ \eta_1\leq n-1}}
\frac{\delta_{\eta}K_{\nu'\eta}(q)}{(q)_{m(\eta)}}.
\end{equation*}
Writing $\nu=(a^n)+\lambda$ with $\lambda$ 
a partition such that $l(\lambda)\leq n-1$ this becomes
\begin{align*}
\gamma_k&=\sum_{a=0}^{\infty}
\sum_{\substack{\lambda\in\Pa\\ l(\lambda)\leq n-1 \\[0.3mm]
|\lambda|\geqq\ell+nk_1\pmd{n}}} 
K_{(a^n)+\lambda,(a^n)+\frac{|\lambda|}{n}(1^n)
-k-\ell\bar{\Lambda}_{n-1}}\\[-3mm]
& \qquad \qquad \qquad \qquad \qquad \qquad \qquad \times
\sum_{\substack{\eta\vdash an+|\lambda| \\ \eta_1\leq n-1}}
\frac{\delta_{\eta}K_{(n^a)\cup\lambda',\eta}(q)}{(q)_{m(\eta)}}  \\
&=\sum_{\substack{\lambda\in\Pa\\ l(\lambda)\leq n-1 \\[0.3mm] 
|\lambda|\geqq\ell+nk_1\pmd{n}}} 
K_{\lambda,\frac{|\lambda|}{n}(1^n)-k-\ell\bar{\Lambda}_{n-1}} 
\sum_{a=0}^{\infty}\:
\sum_{\substack{\eta\vdash an+|\lambda| \\ \eta_1\leq n-1}}
\frac{\delta_{\eta}K_{(n^a)\cup\lambda',\eta}(q)}{(q)_{m(\eta)}},
\end{align*}
where the second equality follows from \eqref{square}.
Interchanging the sums over $\eta$ and $a$ it follows that $a$ 
is fixed by $a=(|\eta|-|\lambda|)/n$, resulting in the claim of the lemma.
\end{proof}

\begin{theorem}\label{thm2}
For $\mu\in P_{+}$ and $\nu\in\Pa$ such that $\nu_1\leq N$ and 
$|\nu|\bar{\Lambda}_1-\mu\in Q$,
\begin{equation*}
\sum_{\substack{\lambda\in\Pa\\l(\lambda)\leq n-1\\[0.3mm]
|\lambda|\bar{\Lambda}_1-\mu\in Q}}
K_{\lambda,\frac{|\lambda|}{n}(1^n)+\mu}
\sum_{\substack{\eta\in\Pa\\ \eta_1\leq n-1 \\[0.3mm]
|\eta|\geqq|\lambda|\pmd{n}}}
\frac{F_{\eta\nu}(q)K_{(n^{(|\eta|-|\lambda|)/n})\cup \lambda',\eta}(q)}
{(q)_{m(\eta)}}=C_{\mu\nu}(q).
\end{equation*}
\end{theorem}
Before proving this, let us show how combined with Lemma~\ref{lem1}
it implies Theorem~\ref{thm1}.
By \eqref{Ksymm}, $\mu$ on the left may be replaced by $\pi(\mu)$.
After this change we choose
$\mu=\ell\bar{\Lambda}_1-\pi(k)$ with $k\in Q\cap P_{+}$
and $\ell\in\Z_{+}$, and use that
$\pi(\bar{\Lambda}_1)=-\bar{\Lambda}_{n-1}$ and 
$(\bar{\Lambda}_1)_a-(\bar{\Lambda}_1)_{a+1}=\delta_{a,1}$.
Hence
\begin{multline}\label{corK}
\sum_{\substack{\lambda\in\Pa\\l(\lambda)\leq n-1\\[0.3mm] 
|\lambda|\geqq\ell+nk_1\pmd{n}}}
K_{\lambda,\frac{|\lambda|}{n}(1^n)-k-\ell\bar{\Lambda}_{n-1}}\\[-2mm]
\times
\sum_{\substack{\eta\in\Pa\\ \eta_1\leq n-1 \\[0.3mm] 
|\eta|\geqq|\lambda|\pmd{n}}}
\frac{F_{\eta\nu}(q)K_{(n^{(|\eta|-|\lambda|)/n})\cup \lambda',\eta}(q)}
{(q)_{m(\eta)}}=C_{\ell\bar{\Lambda}_1-\pi(k),\nu}(q)
\end{multline}
with $|\nu|\equiv\ell\pmod{n}$.
Comparing this with Lemma~\ref{lem1} yields Theorem~\ref{thm1}.

The proof of Theorem~\ref{thm2} requires the following two propositions.
\begin{proposition}\label{prop1}
For $\lambda,\nu\in\Pa$ such that $l(\lambda)\leq n$, $\nu_1\leq N$
and $|\lambda|=|\nu|$,
\begin{equation*}
\lim_{M\to\infty} q^{-nN\binom{M}{2}-M|\nu|}
K_{NM(1^n)+\lambda,(N^{nM})\cup\nu}(q)
=\sum_{\substack{\eta\in\Pa\\ \eta_1\leq n-1\\[0.3mm]|
\eta|\equiv|\lambda|\pmd{n}}} 
\frac{F_{\eta\nu}(q)K_{\xi'\eta}(q)}
{(q)_{m(\eta)}},
\end{equation*}
where $\xi=\frac{|\eta|-|\lambda|}{n}(1^n)+\lambda$.
\end{proposition}
\begin{proposition}\label{prop2}
For $\lambda,\nu\in\Pa$ such that $l(\lambda)\leq n$, $\nu_1\leq N$
and $|\lambda|=|\nu|$,
\begin{equation*}
\lim_{M\to\infty} q^{-nN\binom{M}{2}-M|\nu|}
\Su_{NM(1^n)+\lambda,(N^{nM})\cup\nu}(q)=C_{\mu\nu}(q),
\end{equation*}
where $\mu=\sum_{a=1}^{n-1}(\lambda_a-\lambda_{a+1})\bar{\Lambda}_a\in P_{+}$.
\end{proposition}
Proposition~\ref{prop1} is \cite[Prop. 6.4]{HKKOTY99} and
Proposition~\ref{prop2} follows from \cite[Props. 4.3, 4.6 and 5.8]{HKKOTY99}.
Before giving proofs we first show how these results imply Theorem~\ref{thm2}.
\begin{proof}[Proof of Theorem~\ref{thm2}]
For $\mu\in P$ and $\nu\in\Pa$ such that $\nu_1\leq N$ and 
$|\nu|\bar{\Lambda}_1-\mu\in Q$, consider the expression
\begin{multline}\label{I}
I_{\mu\nu}(q)=\lim_{M\to\infty}q^{-nN\binom{M}{2}-M|\nu|} \\[2mm] \times
\sum_{\substack{\lambda\in\Pa\\l(\lambda)\leq n-1\\[0.3mm]
|\lambda|\bar{\Lambda}_1-\mu\in Q}}
K_{\lambda,\frac{|\lambda|}{n}(1^n)+\mu} 
K_{\frac{nNM+|\nu|-|\lambda|}{n}(1^n)+\lambda,(N^{nM})\cup\nu}(q).
\end{multline}
By Proposition~\ref{prop1} with $\lambda\to \frac{1}{n}(|\nu|-|\lambda|)
(1^n)+\lambda$ it is readily found that
\begin{equation*}
I_{\mu\nu}(q)=
\sum_{\substack{\lambda\in\Pa\\l(\lambda)\leq n-1\\[0.3mm]
|\lambda|\bar{\Lambda}_1-\mu\in Q}}
K_{\lambda,\frac{|\lambda|}{n}(1^n)+\mu} 
\sum_{\substack{\eta\in\Pa\\ \eta_1\leq n-1 \\[0.3mm]
|\eta|\geqq|\lambda|\pmd{n}}}
\frac{F_{\eta\nu}(q)K_{(n^{(|\eta|-|\lambda|)/n})\cup \lambda',\eta}(q)}
{(q)_{m(\eta)}}.
\end{equation*}
Here we have used that $\xi=\frac{|\eta|-|\lambda|}{n}(1^n)+\lambda$
with $l(\lambda)\leq n-1$ is a partition iff $|\eta|\geqq|\lambda|\pmod{n}$.
Hence we may assume this condition and write $\xi'= 
(n^{(|\eta|-|\lambda|)/n})\cup \lambda'$.

In order to prove Theorem~\ref{thm2} it remains to be shown that 
$I_{\mu\nu}=C_{\mu\nu}$ if $\mu\in P_{+}$.
To this end we use the properties \eqref{square} and \eqref{Kzero}
of the Kostka numbers to identify the summand of \eqref{I}
with a completely symmetric supernomial;
\begin{align*}
\sum_{\substack{\lambda\in\Pa\\l(\lambda)\leq n-1\\[0.3mm]
|\lambda|\bar{\Lambda}_1-\mu\in Q}}&
K_{\lambda,\frac{|\lambda|}{n}(1^n)+\mu}
K_{\frac{nNM+|\nu|-|\lambda|}{n}(1^n)+\lambda,(N^{nM})\cup\nu}(q) \\
&=\sum_{\substack{\lambda\in\Pa\\l(\lambda)\leq n-1\\[0.3mm]
|\lambda|\bar{\Lambda}_1-\mu\in Q}}
K_{\frac{nNM+|\nu|-|\lambda|}{n}(1^n)+\lambda,\frac{nNM+|\nu|}{n}(1^n)+\mu} 
\\[-4mm] & \hspace{3.5cm} \times
K_{\frac{nNM+|\nu|-|\lambda|}{n}(1^n)+\lambda,(N^{nM})\cup\nu}(q) \\[3mm]
&=\sum_{\substack{\eta\vdash |\nu|+nNM \\ l(\eta)\leq n}}
K_{\eta,\frac{nNM+|\nu|}{n}(1^n)+\mu}K_{\eta,(N^{nM})\cup\nu}(q) \\
&=\Su_{\frac{nNM+|\nu|}{n}(1^n)+\mu,(N^{nM})\cup\nu}(q).
\end{align*}
The last equality follows from \eqref{SKK} and the fact that the
restriction $l(\lambda)\leq n$ in the second-last line may be dropped 
thanks to \eqref{Ksymm} and \eqref{Kzero}.

By Proposition~\ref{prop2} with $\lambda\to\frac{|\nu|}{n}(1^n)+\mu$
(which is a partition iff $\mu\in P_{+}$) it thus follows that for 
$\mu\in P_{+}$ there holds $I_{\mu\nu}=C_{\mu\nu}$.
\end{proof}

\subsection{Proof of Proposition \ref{prop1}}\label{sec33}

As mentioned previously, Proposition~\ref{prop1} is due to 
Hatayama \textit{et al.} \cite[Prop. 6.4]{HKKOTY99},
who we will closely follow in our proof.
The only significant difference is that in \cite{HKKOTY99} only the 
case $\nu=\emptyset$ is treated in detail.

Key to Proposition~\ref{prop1} are two so-called fermionic representations 
of the Kostka polynomials. The first of these is due to Kirillov and 
Reshetikhin \cite[Thm. 4.2]{KR88}.
Let $T_{j,k}=\min\{j,k\}$.
\begin{proposition}\label{prop3}
For $\lambda,\mu\in\Pa$ such that $|\lambda|=|\mu|$ and $l(\lambda)\leq n$,
\begin{multline}\label{Kferm}
K_{\lambda\mu}(q)=q^{n(\mu)}\sum 
q^{\frac{1}{2}\sum_{a,b=1}^{n-1}\sum_{j,k\geq 1}m_j^{(a)}
C_{a,b}T_{j,k}m_k^{(b)}}\\ \times 
q^{-\sum_{i=1}^{l(\mu)}\sum_{j\geq 1} T_{\mu_i,j} m_j^{(1)}}
\prod_{a=1}^{n-1}\prod_{j\geq 1}
\qbinsb{m_j^{(a)}+p_j^{(a)}}{m_j^{(a)}}_q.
\end{multline}
Here $p_j^{(a)}$ is given by 
\begin{equation}\label{p2}
p_j^{(a)}=\delta_{a,1}\sum_{i=1}^{l(\mu)}T_{j,\mu_i}
-\sum_{b=1}^{n-1}\sum_{k\geq 1}C_{a,b}T_{j,k}m_k^{(b)}
\end{equation}
and the sum is over $m_j^{(a)}\in\Z_{+}$ for $a\in [n-1]$ and $j\geq 1$, 
subject to
\begin{equation}\label{rest2}
\sum_{j\geq 1}j m_j^{(a)}=\sum_{b=a+1}^n \lambda_b.
\end{equation}
\end{proposition}

Using (i) the duality \eqref{Kdual}, (ii) the analogue of 
Proposition~\ref{prop3} for generalized Kostka polynomials, conjectured 
in \cite[Conj. 6]{KS02} and \cite[Conj. 8.3]{SW99} and proven in 
\cite[Thm 2.10]{KSS02}, and (iii) \eqref{KtoK}, a second fermionic 
representation of the Kostka polynomials arises~\cite[Prop. 5.6]{HKKOTY99}.
\begin{proposition}\label{prop4}
For $\lambda,\eta\in\Pa$ such that $|\lambda|=|\eta|$, $l(\lambda)\leq n$ 
and $\eta_1\leq n-1$
\begin{equation}\label{Kd}
K_{\lambda'\eta}(q)=\sum q^{\frac{1}{2}\sum_{a,b=1}^{n-1}
\sum_{j,k\geq 1}m_j^{(a)}C_{a,b}T_{j,k}m_k^{(b)}}
\prod_{a=1}^{n-1}\prod_{j\geq 1}
\qbinsb{m_j^{(a)}+p_j^{(a)}}{m_j^{(a)}}_q.
\end{equation}
Here $p_j^{(a)}$ is given by 
\begin{equation}\label{p3}
p_j^{(a)}=m_a(\eta)-\sum_{b=1}^{n-1}\sum_{k\geq 1}C_{a,b}T_{j,k}m_k^{(b)},
\end{equation}
the sum is over $m_j^{(a)}\in\Z_{+}$ for $a\in [n-1]$ and $j\geq 1$
such that 
\begin{equation}\label{rest3}
\sum_{j\geq 1}j m_j^{(a)}
=\ip{\bar{\Lambda}_a}{\mu}
-\sum_{b=1}^a \Bigl(\lambda_b-\frac{|\lambda|}{n}\Bigr)
\end{equation}
and $\mu=\sum_{a=1}^{n-1}m_a(\eta)\bar{\Lambda}_a$.
\end{proposition}
As an example, let us calculate $K_{(4,3,1),(3,3,2)}(q)$.
Taking the Kirillov--Reshetikhin representation we have to 
compute all solutions to $\sum_{j\geq 1}jm_j^{(1)}=4$ and
$\sum_{j\geq 1}m_j^{(2)}=1$. The first of these equations has
five solutions corresponding to the five partitions of $4$.
The second equation has the unique solution $m_j^{(2)}=\delta_{j,1}$.
Calculating the corresponding $p_j^{(a)}$ for each of the five solutions
using \eqref{p2}, it turns out that only one of the five
has all $p_j^{(a)}$ nonnegative.
Hence only this solution, given by $m_1^{(1)}=m_3^{(1)}=m_1^{(2)}=1$,
$m_j^{(a)}=0$ otherwise and $p_j^{(1)}=\chi(j\geq 2)$, $p_j^{(2)}=p_j^{(1)}+
\chi(j\geq 3)$ contributes to the sum
yielding $K_{(4,3,1),(3,3,2)}(q)=q\qbins{2}{1}_q=q+q^2$.

Taking the representation of Proposition~\ref{prop4} we have 
$\lambda'=(4,3,1)$, $\eta=(3,3,2)$ and hence $\lambda=(3,2,2,1)$, 
$m(\eta)=(0,1,2)$. We may thus take $n=4$ yielding $\mu=(1,1,0,-2)$
and leading to the equations $\sum_{j\geq 1}j m_j^{(1)}=0$ and
$\sum_{j\geq 1}j m_j^{(a)}=1$ for $a\in\{2,3\}$.
This has the unique solution $m_1^{(2)}=m_1^{(3)}=1$ and
$m_j^{(a)}=0$ otherwise. From \eqref{p3} it then follows that
$p_j^{(1)}=p_j^{(3)}=1$ and $p_j^{(2)}=0$ for all $j\geq 1$,
so that, once more, $K_{(4,3,1),(3,3,2)}(q)=q\qbins{2}{1}_q=q+q^2$.

To now prove Proposition \ref{prop1}, take the Kostka polynomial 
in the representation of Proposition~\ref{prop3}, assume that $\mu_1\leq N$,
and replace 
\begin{equation}\label{rep}
\lambda\to NM(1^n)+\lambda \quad \text{and} \quad \mu\to (N^{nM})\cup\nu
\end{equation}
with $\lambda,\nu\in\Pa$ such that $l(\lambda)\leq n$, $\nu_1\leq N$ and 
$|\lambda|=|\nu|$.
By this change equations \eqref{p2} and \eqref{rest2} become
\begin{equation}\label{nshift}
p_j^{(a)}=\delta_{a,1}\sum_{i=1}^{l(\nu)}T_{j,\nu_i}
-\sum_{b=1}^{n-1}\sum_{k\geq 1}C_{a,b}T_{j,k}
\bigl(m_k^{(b)}-v_k^{(b)}\bigr)
\end{equation}
and
\begin{equation}\label{rest4}
\sum_{j\geq 1}j \bigl(m_j^{(a)}-v_j^{(a)}\bigr)
=\sum_{b=a+1}^n \lambda_b,
\end{equation}
where $v_j^{(a)}=M(n-a)\delta_{j,N}$. 
Without loss of generality we may assume that $p_N^{(a)}\geq 0$,
and we define $\eta\in\Pa$ such that $\eta_1\leq n-1$ by 
$m_a(\eta)=p_N^{(a)}$ and set 
$\mu=\sum_{a=1}^{n-1}m_a(\eta)\bar{\Lambda}_a\in P_{+}$.
Then \eqref{nshift} implies
\begin{equation}\label{zeta}
m_a(\eta)=|\nu|\delta_{a,1}
-\sum_{b=1}^{n-1}\sum_{j\geq 1}C_{a,b}T_{N,j}
\bigl(m_j^{(b)}-v_j^{(b)}\bigr).
\end{equation}
By $|\eta|=\sum_{a=1}^{n-1}am_a(\eta)$ and 
$\sum_{a=1}^{n-1}aC_{a,b}=n\delta_{b,n-1}$ it follows that 
$\eta\equiv|\nu|\pmod{n}$. The equations \eqref{nshift}, \eqref{zeta} and 
$T_{j,k}=\chi(k\in[N-1])\bar{C}^{-1}_{j,k}+jT_{N,k}/N$ 
(true for $k\geq 0$ and $0\leq j\leq N$) further yield 
\begin{align}\label{pja1}
p_j^{(a)} 
&=\delta_{a,1}\sum_{i=1}^{l(\nu)}T_{j,\nu_i}
-\sum_{b=1}^{n-1}C_{a,b}\biggl[\:\sum_{k=1}^{N-1}\bar{C}^{-1}_{j,k}m_k^{(b)}
+\frac{j}{N}\sum_{k\geq 1}T^{-1}_{N,k}
\bigl(m_k^{(b)}-v_k^{(b)}\bigr)\biggr] \\
&=\delta_{a,1}\sum_{i=1}^{l(\nu)}\bar{C}^{-1}_{j,\nu_i}
+\frac{j}{N}m_a(\eta)
-\sum_{b=1}^{n-1}\sum_{k=1}^{N-1}C_{a,b}\bar{C}^{-1}_{j,k}m_k^{(b)}\notag
\end{align}
for $j\in[N-1]$. 
Similarly, by \eqref{nshift}, $\nu_1\leq N$ and
$T_{j+N,k}=\chi(k>N)T_{j,k-N}+T_{k,N}$ (true for $j\geq 0$) one finds
\begin{align}\label{pja2}
p_{j+N}^{(a)}
&=|\nu|\delta_{a,1}-\sum_{b=1}^{n-1}C_{a,b}\biggl[\:
\sum_{k\geq 1}T_{j,k}m_{k+N}^{(b)}
+\sum_{k\geq 1}T_{k,N}\bigl(m_k^{(b)}-v_k^{(b)}\bigr)\biggr] \\
&=m_a(\eta)-\sum_{b=1}^{n-1}\sum_{k\geq 1}C_{a,b}T_{j,k}m_{k+N}^{(b)}
\notag
\end{align}
for $j\geq 1$. Note that \eqref{pja1} for $j=N$ and \eqref{pja2} for
$j=0$ correspond to the tautology $p_N^{(a)}=m_a(\eta)$.

Equation \eqref{zeta}, the definition of $\mu$ and $|\nu|=|\lambda|$
may also be applied to rewrite the restriction \eqref{rest4}. 
Namely, if in \eqref{zeta} we replace $a\to d$, then multiply by 
$C^{-1}_{a,d}$ and sum over $d$, and finally subtract \eqref{rest4}
from the resulting equation, we find
\begin{equation}\label{zf}
\sum_{j\geq 1}jm_{j+N}^{(a)}=\ip{\bar{\Lambda}_a}{\mu}-
\sum_{b=1}^a\Bigl(\lambda_b-\frac{|\lambda|}{n}\Bigr).
\end{equation}

Finally consider the exponent $E$ of $q$ on the right of \eqref{Kferm}
after the replacements \eqref{rep}. From \eqref{zeta} and \eqref{nab}
it follows that
\begin{align*}
E&=nN\binom{M}{2}+M|\nu|+n(\nu)-\frac{(n-1)|\nu|^2}{2nN}
+\frac{\ip{\mu}{\mu}}{2N} \\
& \qquad +\frac{1}{2}\sum_{a,b=1}^{n-1}\sum_{j,k=1}^{N-1}
m_j^{(a)}C_{a,b}\bar{C}^{-1}_{j,k} \Bigl(m_k^{(b)}-
\sum_{i=1}^{l(\nu)}\delta_{a,1}\delta_{b,1}\delta_{k,\nu_i}\Bigr) \\
& \qquad +\frac{1}{2}\sum_{a,b=1}^{n-1}\sum_{j,k\geq 1}
C_{a,b}\bar{C}^{-1}_{j,k}m_{j+N}^{(a)}m_{k+N}^{(b)}.
\end{align*}
This result does not contain $m_N^{(a)}$,
and in order to take the large $M$ limit we replace 
$m_j^{(a)}\to m_j^{(a)}+v_j^{(a)}=m_j^{(a)}+M(n-a)\delta_{j,N}$.
Then the $M$-dependence of \eqref{zeta} drops out
and the only occurrence of $M$ in the above equations is in the expression 
for $E$, which contains the term $nN\binom{M}{2}+M|\nu|$.
We now define $\bar{m}$ by \eqref{m} with $m_j^{(a)}$ therein
replaced by $m_{N-j}^{(a)}$ (i.e., $\bar{m}_j^{(a)}=m_{N-j}^{(a)}$),
and define $\tilde{m}_j^{(a)}=m_{N+j}^{(a)}$. 
Correspondingly, define $\bar{p}$ ($\tilde{p}_j^{(a)}$)
by \eqref{p1} (\eqref{p3}) with $m$ replaced by $\bar{m}$ 
($\tilde{m}$).
Then \eqref{pja1} (\eqref{pja2}) is nothing but \eqref{p1} (\eqref{p3})
with $m$ and $p$ replaced by $\bar{m}$ and $\bar{p}$ 
($\tilde{m}$ and $\tilde{p}$), and \eqref{zeta} becomes (after also
multiplying by $\bar{\Lambda}_a$ and summing over $a$)
\begin{equation}\label{mNres}
\sum_{a=1}^{n-1}\sum_{j=1}^{N-1}j\bar{m}_j^{(a)}\alpha_a=
\mu-|\nu|\bar{\Lambda}_1
+N\sum_{a=1}^{n-1}\alpha_a\biggl[m_N^{(a)}+\sum_{j=1}^{N-1}\bar{m}_j^{(a)}
+\sum_{j\geq 1}\tilde{m}_j^{(a)}\biggr].
\end{equation}
Now writing the expressions \eqref{F} and \eqref{Kd} without the 
summations on the left as $F_{\eta\nu;m}(q)$ and
$K_{\lambda'\eta;\{m_j^{(a)}\}}(q)$ results in
\begin{multline}\label{KmFm}
q^{-nN\binom{M}{2}-M|\nu|}
K_{NM(1^n)+\lambda,(N^{nM})\cup\nu}(q) \\
=\sum F_{\eta\nu;\bar{m}}(q)
K_{\xi'\eta;\{\tilde{m}_j^{(a)}\}}(q)
\prod_{a=1}^{n-1}\qbinsb{m_N^{(a)}+m_a(\eta)+(n-a)M}{m_a(\eta)}_q,
\end{multline}
where the sum is over $\bar{m}$, the $\tilde{m}_j^{(a)}$ and $m_N^{(a)}$
subject to the following restrictions: (i)  
$\mu$ is fixed by \eqref{zf}, i.e., by \eqref{rest3}
with $m$ replaced by $\tilde{m}$ but must be an element of $P_{+}$ (so that
$m_a(\eta)=\ip{\mu}{\alpha_a}\geq 0$, and (ii) \eqref{mNres} must holds.
To understand the occurrence of $\xi$ in the above, note that it is
the unique partition such that $|\xi|=|\eta|$ and
$\xi_b-|\xi|/n=\lambda_b-|\lambda|/n$.

To complete the proof, replace the above sum by the equivalent sum over 
$\mu\in P_{+}$, $\bar{m}$ and the $\tilde{m}_j^{(a)}$ where the latter
are subject to the restriction \eqref{rest3} with 
$m\to \tilde{m}$, and $\bar{m}$ is restricted by the condition
that \eqref{mNres} must fix the $m_N^{(a)}$ to be integers. 
If we now take the large $M$ limit, the product on the right side of
\eqref{KmFm} yields $1/(q)_{m(\eta)}$ 
eliminating all $m_N^{(a)}$-dependence
of the summand. Hence  we may replace \eqref{mNres}
by \eqref{rest1} with $m\to\bar{m}$. Finally replacing the sum over
$\mu$ by a sum over $\eta$ yields Proposition~\ref{prop1}.

\subsection{Proof of Proposition \ref{prop2}}\label{sec34}
This proposition follows from \cite[Props. 4.3, 4.6 and 5.8]{HKKOTY99}
and our proof below does not significantly differ from the one given
by Hatayama \textit{et al.}

It was shown in \cite[Prop. 5.1]{HKKOTY99} that for $|\lambda|=|\mu|$
and $\lambda\in\Z_{+}^n$ the completely symmetric supernomial
can be expressed as
\begin{equation*}
\Su_{\lambda\mu}(q)=
\sum_{\nu} q^{\sum_{a=0}^{n-1}\sum_{j=1}^{\mu_1}
\Bigl(\genfrac{}{}{0pt}{2}{\nu_j^{(a+1)}-\nu_j^{(a)}}{2}\Bigr)}
\prod_{a=1}^{n-1}\prod_{j=1}^{\mu_1}
\qbinsb{\nu_j^{(a+1)}-\nu_{j+1}^{(a)}}{\nu_j^{(a+1)}-\nu_j^{(a)}}_q,
\end{equation*}
where the sum is over sequences of partitions 
$\nu=(\nu^{(1)},\dots,\nu^{(n-1)})$ such that
$\emptyset=\nu^{(0)}\subset\nu^{(1)}\subset\dots\subset\nu^{(n)}=\mu'$
and $|\nu^{(a+1)}|-|\nu^{(a)}|=\lambda_{n-a}$.
For $n=2$ this is equivalent to \cite[Eqs. (2.9)--(2.10)]{SW98b}.

Now assume that $\mu_1=N$ and introduce the variables
$m_j^{(a)}=\nu_{N-j}^{(n-a)}-\nu_{N-j+1}^{(n-a)}$ for $j\in\{0,\dots,N-1\}$.
Also defining the vector $m$ as in \eqref{m} gives
\begin{multline*}
\Su_{\lambda\mu}(q)=q^{n(\mu)-\frac{(n-1)|\mu|^2}{2nN}
+\frac{1}{2N}\sum_{a=1}^n \bigl(\lambda_a^2-\frac{|\lambda|^2}{n^2}\bigr)} \\
\times\sum q^{\frac{1}{2}m (C\otimes\bar{C}^{-1})m
-\sum_{i=1}^{l(\mu)}(e_1\otimes \bar{e}_{N-\mu_i})(I\otimes\bar{C}^{-1})m}
\prod_{a=1}^{n-1}\prod_{j=0}^{N-1}
\qbinsb{m_j^{(a)}+p_j^{(a)}}{m_j^{(a)}}_q,
\end{multline*}
where the sum is over $m$ and the $m_0^{(a)}$
such that $\sum_{j=0}^{N-1}(N-j)m_j^{(a)}=\sum_{b=a+1}^n\lambda_b$.
The auxiliary variables $p_j^{(a)}$ are defined as
$p_j^{(a)}=\sum_{k=0}^j (m_k^{(a-1)}-m_k^{(a)})+\delta_{a,1}\mu'_{N-j}$
(with $m_j^{(0)}:=0$).
Making the replacements \eqref{rep} in the above result, one finds
\begin{multline*}
q^{-nN\binom{M}{2}-M|\nu|}
\Su_{NM(1^n)+\lambda,(N^{nM})\cup\nu}(q)
=q^{n(\nu)-\frac{(n-1)|\nu|^2}{2nN}+\frac{\ip{\mu}{\mu}}{2N}} \\
\times \sum q^{\frac{1}{2}m (C\otimes\bar{C}^{-1})m
-\sum_{i=1}^{l(\nu)}(e_1\otimes \bar{e}_{N-\nu_i})
(I\otimes\bar{C}^{-1})m}
\prod_{a=1}^{n-1}\prod_{j=0}^{N-1}
\qbinsb{m_j^{(a)}+p_j^{(a)}}{m_j^{(a)}}_q.
\end{multline*}
Here $\lambda,\nu$ and $\mu$ are as in Proposition~\ref{prop2},
the sum is over $m$ and the $m_0^{(a)}$ such that
\begin{equation}\label{rest5}
\sum_{a=1}^{n-1}\sum_{j=1}^{N-1}jm_j^{(a)}\alpha_a=
\mu-|\nu|\bar{\Lambda}_1+N\sum_{a=1}^{n-1}\alpha_a
\biggl[\sum_{j=0}^{N-1}m_j^{(a)}-(n-a)M\biggr]
\end{equation}
and $p_j^{(a)}=\sum_{k=0}^j (m_k^{(a-1)}-m_k^{(a)})+
\delta_{a,1}(nM+\nu'_{N-j})$.
To arrive at the exponent of $q$ on the right we have used
\eqref{nab} and $\sum_{a=1}^n(\lambda_a^2-|\lambda|^2/n^2)=
\ip{\mu}{\mu}$, and to obtain the restriction in the sum we have used that
$\sum_{a=1}^{n-1}\alpha_a\sum_{b=a+1}^n\lambda_b=|\nu|\bar{\Lambda}_1-\mu$.
Since the exponent of $q$ does not contain $M$, $m_0^{(a)}$ and $p_j^{(a)}$
it is straightforward to let $M$ tend to infinity.
All we need to do is replace $m_0^{(a)}\to (n-a)M+m_0^{(a)}$.
Then the product over the $q$-binomials becomes
\begin{equation*}
P:=\prod_{a=1}^{n-1}\biggl(
\qbinsb{(n-a)M+m_0^{(a)}+p_0^{(a)}}{(n-a)M+m_0^{(a)}}_q
\prod_{j=1}^{N-1}
\qbinsb{m_j^{(a)}+p_j^{(a)}}{m_j^{(a)}}_q \:\biggr),
\end{equation*}
with $p_j^{(a)}=\sum_{k=0}^j (m_k^{(a-1)}-m_k^{(a)})+M+
\delta_{a,1}\nu'_{N-j}$ so that 
$\lim_{M\to\infty}P=1/[(q)_{\infty}^{n-1}(q)_m]$ which is independent of 
the $m_0^{(a)}$.
Moreover, the term $(n-a)M$ on the right of \eqref{rest5} is cancelled
so that summing over $m$ and the $m_0^{(a)}$
is equivalent to summing over $m$ subject to the restriction \eqref{rest1}.

\section{Bosonic representation of Theorem~\ref{thm1}}\label{sec4}
The aim of this section is to provide an alternative representation
(for a special case) of Theorem~\ref{thm1}, in which $\gamma_k$ and 
$\delta_{\eta}$ are expressed as an A$_{n-1}^{(1)}$ string function and an 
A$_{n-1}^{(1)}$ configuration sum, both at level $N$.
The advantage of this ``bosonic'' representation of the theorem
is that it suggests a fractional-level generalization in which
the integer $N$ is replaced by a rational number.

Let $\Lambda_0,\dots,\Lambda_{n-1}$ be the fundamental weights of
the affine Lie algebra A$_{n-1}^{(1)}$ \cite{Kac90}, 
and define $\Lambda_i$ for
all integers $i$ by $\Lambda_i=\Lambda_j$ if $i\equiv j\pmod{n}$.
The set of level-$N$ dominant integral weights $P_{+}^{N}$ is the
set of weights of the form $\Lambda=\sum_{i=0}^{n-1}a_i\Lambda_i$ with
$a_i\in\Z_{+}$ and $\sum_{i=0}^{n-1}a_i=N$.
The classical part of a weight $\Lambda\in P_{+}^N$ will be 
denoted by $\bar{\Lambda}$, consistent with our definition 
of the fundamental weights of A$_{n-1}$. Hence
$\bar{\Lambda}_0=\emptyset$.
When $N=0$, $P_{+}^N=\emptyset$, but it will be convenient
to follow the convention of partitions, and to somewhat ambiguously
write $P_{+}^0=\{\emptyset\}$, with $\emptyset$ the ``empty'' level-$0$
weight.
The modular anomaly $m_{\Lambda}$ for $\Lambda\in P_{+}^N$ is given by
\begin{equation*}
m_{\Lambda}=\frac{||\bar{\Lambda}+\rho||^2}{2(n+N)}-\frac{||\rho||^2}{2n},
\end{equation*}
with $||\rho||^2=n\dim(\text{A}_{n-1})/12=n(n^2-1)/12$.
We will also employ a shifted modular anomaly
\begin{equation*}
\bar{m}_{\Lambda}=\frac{||\bar{\Lambda}+\rho||^2}{2(n+N)}.
\end{equation*}

We now define a level-$N$, A$_{n-1}^{(1)}$ configuration sum 
$X_{\eta,\Lambda,\Lambda'}$ by
\begin{equation}\label{XN}
X_{\eta,\Lambda,\Lambda'}(q)=q^{-\bar{m}_{\Lambda}}
\sum_{\sigma\in S_n}\epsilon(\sigma)
\sum_{\substack{k\in (n+N)Q+\sigma(\bar{\Lambda}+\rho)}}
q^{\frac{\ip{k}{k}}{2(n+N)}}
\qbin{m(\eta)}{\frac{|\eta|}{n}(1^n)+k-\Lambda'-\rho}
\end{equation}
for $\Lambda,\Lambda'\in P_{+}^N$ and $\eta\in\Pa$ such that $\eta_1\leq n-1$
and $|\eta|\bar{\Lambda}_1-\bar{\Lambda}+\bar{\Lambda}'\in Q$.

Using a sign-reversing involution, Schilling and Shimozono \cite{SS00}
showed that $X_{\eta,\Lambda,N\Lambda_0}(q)$ is a special case of the 
level-$N$ restricted generalized Kostka polynomial. Equating this with 
the fermionic representation for this polynomial \cite{SS01} leads to 
identities for $X_{\eta,\Lambda,N\Lambda_0}$, which 
for a special subset of $\Lambda\in P_{+}^N$ can be stated 
as follows \cite{SS00,SS01} (see also \cite{OSS01}).
\begin{theorem}\label{thm3}
Let $\Lambda=(N-s)\Lambda_r+s\Lambda_{r+1}$
with $r\in\{0,\dots,n-1\}$ and $s\in\{0,\dots,N\}$ and
let $\eta\in\Pa$ such that $\eta_1\leq n-1$ and
$|\eta|\bar{\Lambda}_1-\bar{\Lambda}\in Q$. Then
\begin{equation}\label{SS}
X_{\eta,\Lambda,N\Lambda_0}(q)=F_{\eta,(sN^r)}(q).
\end{equation}
\end{theorem}
Note that the second restriction on $\eta$ is equivalent to
$|\eta|\equiv rN+s\pmod{n}$.
Also note that from \eqref{FN1}, \eqref{FCN0} and \eqref{SS}
follow the particularly simple (to state, \textit{not} to prove)
identities \cite[Eq. (6.5)]{SS00} 
\begin{equation}\label{X1}
X_{\eta,\Lambda_i,\Lambda_0}(q)=
q^{\frac{1}{2}\ip{\mu}{\mu}-\frac{1}{2}||\bar{\Lambda}_i||^2}
\end{equation}
for $|\eta|\equiv i\pmod{n}$ with 
$\mu=\sum_{a=1}^{n-1}m_a(\eta)\bar{\Lambda}_a$, and \cite[Eq. (6.6)]{SS00} 
\begin{equation}\label{X0}
X_{\eta,\emptyset,\emptyset}(q)=\delta_{\eta,\emptyset}
\end{equation}
for $|\eta|\equiv 0\pmod{n}$.

Next we shall relate $C_{\mu,(sN^r)}$ to the string functions of 
A$_{n-1}^{(1)}$.
For $\lambda\in P$ and $m$ a positive integer, the classical A$_{n-1}$ theta
function of degree $m$ and characteristics $\lambda$ is defined by
\begin{equation}\label{Tdef}
\Theta_{\lambda,m}(x;q)=\sum_{k\in mQ+\lambda} q^{\frac{1}{2m}\ip{k}{k}}x^k.
\end{equation}
According to the Weyl--Kac formula \cite{Kac74,Kac90} 
the (normalized) character $\chi_{\Lambda}$
of the integrable highest weight module $L(\Lambda)$ of highest weight
$\Lambda\in P_{+}^N$ can be expressed in terms of theta functions as
\begin{equation*}
\chi_{\Lambda}(x;q)=\frac{
\sum_{\sigma\in S_n}\epsilon(\sigma)
\Theta_{\sigma(\bar{\Lambda}+\rho),n+N}(x;q)}
{\sum_{\sigma\in S_n}\epsilon(\sigma)
\Theta_{\sigma(\rho),n}(x;q)}.
\end{equation*}
The level-$N$, A$_{n-1}^{(1)}$ string functions $C_{\mu,\Lambda}$ for
$\mu\in P$ and $\Lambda\in P_{+}^N$ arise as coefficients in the expansion
of $\chi_{\Lambda}$ in terms of degree-$N$ theta functions \cite{Kac80}:
\begin{equation}\label{CN}
\chi_{\Lambda}(x;q)=\sum_{\mu\in P/NQ}
C_{\mu,\Lambda}(q) \Theta_{\mu,N}(x;q).
\end{equation}
Two simple properties of the string functions are \cite{KP84}
\begin{equation*}
C_{\mu,\Lambda}(q)=0 \text{ if $\mu-\bar{\Lambda}\not\in Q$}
\end{equation*}
and
\begin{equation}\label{s1}
C_{\sigma(\mu)+N\alpha,\Lambda}=C_{\mu,\Lambda}
\text{ for $\alpha\in Q$ and $\sigma\in S_n$.}
\end{equation}
For our purposes it will be convenient to also introduce normalized string 
functions $\C_{\mu,\Lambda}$ through
\begin{equation*}
\C_{\mu,\Lambda}(q)=q^{\frac{1}{2N}\ip{\mu}{\mu}-m_{\Lambda}}
C_{\mu,\Lambda}(q).
\end{equation*}

The relevance of the string functions to $C_{\mu\nu}(q)$ of \eqref{C}
lies in the following theorem.
\begin{theorem}\label{thm4}
Let $\Lambda=(N-s)\Lambda_r+s\Lambda_{r+1}$
with $r\in\{0,\dots,n-1\}$ and $s\in\{0,\dots,N\}$ and
let $\mu\in P_{+}$ such that $\mu-\bar{\Lambda}\in Q$. Then
\begin{equation}\label{CC}
\C_{\mu,\Lambda}(q)=C_{\mu,(sN^r)}(q).
\end{equation}
\end{theorem}
This theorem is the culmination of quite a number of results.
For $N=1$ we have according to \eqref{CN1}
\begin{equation}\label{level1}
\C_{\mu,\Lambda_i}(q)=\frac{q^{\frac{1}{2}\ip{\mu}{\mu}-
\frac{1}{2}||\bar{\Lambda}_i||^2}}{(q)_{\infty}^{n-1}}
\end{equation}
for $\mu-\bar{\Lambda}_i\in Q$. This was first obtained by Kac in \cite{Kac78}.
For $n=2$ but arbitrary $N$ the above theorem is due to Lepowsky and 
Primc \cite{LP85} (see also \cite{FS93}).
For general $n$ the theorem was conjectured by Kuniba \textit{et al.} 
\cite{KNS93} for $r=s=0$ and proved by Georgiev \cite{Georgiev95} 
for $r=0$ and $r=n-1$, and by Hatayama \textit{et al.} 
\cite[Prop. 5.8]{HKKOTY99} for general $r$.

If we now take the conjugate Bailey pair of Theorem~\ref{thm1} with 
$\nu=(sN^r)$ and eliminate $F_{\eta\nu}$ and 
$C_{\ell\bar{\Lambda}_1-\pi(k),\nu}$ by 
virtue of Theorems~\ref{thm3} and \ref{thm4}, we arrive at the 
$\Lambda=(N-s)\Lambda_r+s\Lambda_{r+1}$ case of the following conjecture.
\begin{conjecture}\label{conj2}
Let $\ell\in\Z_{+}$, $k\in Q\cap P_{+}$ and $\eta\in\Pa$ such that
$\eta_1\leq n-1$ and $|\eta|\equiv\ell\pmod{n}$.
For $\Lambda\in P_{+}^N$ such that $\ell\bar{\Lambda}_1-\bar{\Lambda}\in Q$ 
\begin{equation*}
\gamma_k=\C_{\ell\bar{\Lambda}_1-\pi(k),\Lambda}(q)\quad\text{and}\quad
\delta_{\eta}=X_{\eta,\Lambda,N\Lambda_0}(q)
\end{equation*}
form an \textup{A}$_{n-1}$ conjugate Bailey pair relative to $q^{\ell}$.
\end{conjecture}
For $\Lambda=(N-s)\Lambda_r+s\Lambda_{r+1}$ 
we can of course claim this as a theorem thanks to 
Theorems~\ref{thm3} and \ref{thm4}.
For $n=2$ a proof of the above result is implicit in \cite{SW98a}. 
In the next section it will be shown that Conjecture~\ref{conj2}
follows by manipulating Corollary~\ref{corinv}.

\section{Fractional-level conjugate Bailey pairs}\label{sec5}
The advantage of Conjecture~\ref{conj2} over Theorem~\ref{thm1} is
that it readily lends itself to further generalization.
In order to describe this we extend our definitions of the
configuration sum $X_{\eta,\Lambda,\Lambda'}$ and string function
$\C_{\mu,\Lambda}$ to fractional levels.

Let $p$ and $p'$ be integers such that $p\geq 1$ and $p'\geq n$,
and fix the level $N$ in terms of $p$ and $p'$ as
\begin{equation}\label{Npp}
N=p'/p-n.
\end{equation}
In contrast with the previous two sections, $N\in\mathbb{Q}$ with $N>-n$.
Given $p$ and $p'$ the set $P^{(p,p')}$ of cardinality $\binom{p'-1}{n-1}$
is defined as the set of weights 
$\Lambda$ of the form $\Lambda=\sum_{i=0}^{n-1}a_i \Lambda_i$ such that 
$\bar{\Lambda}\in P_{+}$ (i.e., $a_i\in\Z_{+}$ for $i\geq 1$),
$\sum_{i=1}^{n-1}a_i\leq p'-n$ and $\sum_{i=0}^{n-1}a_i=N$.
Hence $P_{+}^{(1,p')}=P_{+}^N$, but for $p\geq 2$ the coefficient $a_0$ 
is noninteger if $p$ is not a divisor of $p'$; 
$p'/p-a_0\in\{n,n+1,\dots,p'\}$. Of course, whenever we write $P_{+}^N$
it is assumed that $N\in\Z_{+}$.
Note that the map $\Lambda\to \Lambda+(1-1/p)p'\Lambda_0$ defines a
bijection between $P^{(p,p')}$ and $P_{+}^{p'-n}$.
For $p$ and $p'$ relatively prime
the set $P^{(p,p')}$ is a special subset of the set of
so-called admissible level-$N$ weights \cite{KW88,KW90} in that only 
those weights are included whose classical part is in $P_{+}$. 
We do not know how to include more general admissible weights
in the results of this section.

The definition of the configuration sum \eqref{XN} is most easily 
generalized to fractional $N$, and for 
$\Lambda,\Lambda'\in P^{(p,p')}$ and $\eta\in\Pa$ such that $\eta_1\leq n-1$ 
and $|\eta|\bar{\Lambda}_1-\bar{\Lambda}+\bar{\Lambda}'\in Q$, we define
\begin{equation}\label{Xp}
X_{\eta,\Lambda,\Lambda'}(q)=q^{-\bar{m}_{\Lambda}}
\sum_{\sigma\in S_n}\epsilon(\sigma)
\sum_{\substack{k\in p'Q+\sigma(\bar{\Lambda}+\rho)}}
q^{\frac{\ip{k}{k}}{2(n+N)}}\qbin{m(\eta)}
{\frac{|\eta|}{n}(1^n)+k-\bar{\Lambda}'-\rho}.
\end{equation}
If $|\eta|\bar{\Lambda}_1-\bar{\Lambda}+\bar{\Lambda}'\not\in Q$
we set $X_{\eta,\Lambda,\Lambda'}(q)=0$.

An important remark concerning \eqref{Xp} is in order.
We have not imposed that $p$ and $p'$ be relatively prime,
and as a consequence, $N$ does not uniquely fix $p$ and $p'$.
Hence, given an admissible weight of level $N$, one cannot determine
$p$ and $p'$ . Put differently, the sets $P^{(p,p')}$
and $P^{(t,t')}$ with $(p,p')\neq (t,t')$ are not necessarily
disjoint. For example, $P^{(p,p')}\subseteq P^{(kp,kp')}$ for
$k$ a positive integer.
As a result of this, definition \eqref{Xp} is ambiguous
since the right-hand side depends not only on $N$ (which is
fixed uniquely by $\Lambda$) but also on $p'$.
We however trust that by writing $\Lambda\in P^{(p,p')}$ the reader
will have no trouble interpreting \eqref{Xp}.
Those unwilling to accept the above notation may from now on assume that
$p$ and $p'$ are relatively prime, or should add a superscript $p'$ or even
$(p,p')$ to $X_{\eta,\Lambda,\Lambda'}$.

Before we show how to also generalize \eqref{CN} to yield
fractional-level string functions, we shall list several important 
properties of $X_{\eta,\Lambda,\Lambda'}$, a number of which are conjectural.
\begin{lemma}\label{lem2}
For $\Lambda,\Lambda'\in P^{(p,p')}$,
$X_{\emptyset,\Lambda,\Lambda'}(q)=\delta_{\Lambda,\Lambda'}$.
\end{lemma}
\begin{proof}
{}From \eqref{SKK2} or \eqref{Srev} it follows that 
$\qbins{\emptyset}{\lambda}=\delta_{\lambda,\emptyset}$. Hence
\begin{equation*}
X_{\emptyset,\Lambda,\Lambda'}(q)=q^{\bar{m}_{\Lambda'}-\bar{m}_{\Lambda}}
\sum_{\sigma\in S_n}\epsilon(\sigma)
\sum_{\substack{k\in p'Q+\sigma(\bar{\Lambda}+\rho)}}
\delta_{k,\bar{\Lambda}'+\rho},
\end{equation*}
where we have used that $q^{\ip{k}{k}/(2(n+N))}\delta_{k,\bar{\Lambda}'+\rho}=
q^{\bar{m}_{\Lambda'}}\delta_{k,\bar{\Lambda}'+\rho}$.
The question now is: can we satisfy
$\bar{\Lambda}'+\rho\in p'Q+\sigma(\bar{\Lambda}+\rho)$?
Since $\Lambda,\Lambda'\in P^{(p,p')}$ this simplifies to the problem
of finding for which $\Lambda,\Lambda'$ and $\sigma$ there holds
$\bar{\Lambda}'+\rho=\sigma(\bar{\Lambda}+\rho)$.
Since $\bar{\Lambda}'+\rho\in P_{+}$ and $\sigma(\bar{\Lambda}+\rho)\in P_{+}$
iff $\sigma=(1,2,\dots,n)$, the only solution is given by $\Lambda=\Lambda'$
and $\sigma=(1,2,\dots,n)$, establishing the claim of the lemma.
\end{proof}

For $\Lambda=\sum_{i=0}^{n-1}a_i\Lambda_i\in P^{(p,p')}$ define
$\Lambda^c\in P^{(p,p')}$ as $\Lambda^c=\sum_{i=0}^{n-1}a_i\Lambda_{n-i}$,
and for $\eta=(1^{\zeta_1},2^{\zeta_2},\dots,(n-1)^{\zeta_{n-1}})$
define $\eta^c=(1^{\zeta_{n-1}},2^{\zeta_{n-2}},\dots,(n-1)^{\zeta_1})$.
Note that the partition $\eta^c$ is the complement of $\eta$ with
respect to $(n^{l(\eta)})$, and in our earlier notation of
Section~\ref{sec21}, $\eta^c=\tilde{\eta}_{(n^{l(\eta)})}$. 
\begin{lemma}
For $\Lambda,\Lambda'\in P^{(p,p')}$ and $\eta\in\Pa$ such that
$\eta_1\leq n-1$,
\begin{equation}\label{Xsymm}
X_{\eta,\Lambda,\Lambda'}=X_{\eta^c,\Lambda^c,(\Lambda')^c}.
\end{equation}
\end{lemma}
Observe that the condition
$|\eta|\bar{\Lambda}_1-\bar{\Lambda}+\bar{\Lambda}'\in Q$ implies
that $|\eta^c|\bar{\Lambda}_1-\bar{\Lambda}^c+(\bar{\Lambda}')^c\in Q$ 
as it should.

\begin{proof}
By the symmetry \eqref{Srev} of the antisymmetric supernomials
\begin{equation*}
X_{\eta,\Lambda,\Lambda'}(q)=q^{-\bar{m}_{\Lambda}}
\sum_{\sigma\in S_n}\epsilon(\sigma)
\sum_{\substack{k\in p'Q+\sigma(\bar{\Lambda}+\rho)}}
q^{\frac{\ip{k}{k}}{2(n+N)}}
\qbin{m(\eta^c)}{\frac{|\eta^c|}{n}(1^n)-k+\Lambda'+\rho}.
\end{equation*}
Next use \eqref{Ssymm} with $\sigma=\pi$ and replace $k$ by $-\pi(k)$.
Then
\begin{multline*}
X_{\eta,\Lambda,\Lambda'}(q) \\=q^{-\bar{m}_{\Lambda}}
\sum_{\sigma\in S_n}\epsilon(\sigma)
\sum_{\substack{k\in p'Q-\pi\circ\sigma(\bar{\Lambda}+\rho)}}
q^{\frac{\ip{k}{k}}{2(n+N)}}
\qbin{m(\eta^c)}{\frac{|\eta^c|}{n}(1^n)+k+\pi(\Lambda'+\rho)}.
\end{multline*}
Finally use $\bar{\Lambda}+\rho=-\pi(\bar{\Lambda}^c+\rho)$
(which implies $\bar{m}_{\Lambda}=\bar{m}_{\Lambda^c}$),
$\pi(\bar{\Lambda}'+\rho)=-(\bar{\Lambda}')^c-\rho$
and replace $\sigma$ by $\pi\circ\sigma\circ\pi$ to find \eqref{Xsymm}.
\end{proof}

Based on extensive computer assisted experiments we are led to the
following conjectures.
\begin{conjecture}
Let $p'-p-n+1\geq 0$.
For $\Lambda,\Lambda'\in P^{(p,p')}$ such that
$\Lambda'+(p-1)(p'/p-1)\Lambda_0\in P_{+}^{p'-p-n+1}$, 
$X_{\eta,\Lambda,\Lambda'}(q)$ is a polynomial with nonnegative 
coefficients.
\end{conjecture}
\begin{conjecture}\label{conj3}
For $\Lambda,\Lambda'\in P^{(p,p')}$,
$\lim_{|\eta|\to\infty}X_{\eta,\Lambda,\Lambda'}(q)=0$.
\end{conjecture}
For $n=2$ both conjectures readily follow from results in \cite{SW00a}.
The last conjecture will be crucial later, ensuring the convergence
of certain sums over $X_{\eta,\Lambda,\Lambda'}(q)$.

Next we come to the definition of the fractional-level string functions.
For $\Lambda\in P^{(p,p')}$ define
\begin{equation}\label{KW}
\chi_{\Lambda}(x;q)=
\frac{\sum_{\sigma\in S_n}\epsilon(\sigma)
\Theta_{\sigma(\bar{\Lambda}+\rho),p'}(x;q^p)}
{\sum_{\sigma\in S_n}\epsilon(\sigma)\Theta_{\sigma(\rho),n}(x;q)}.
\end{equation}
For $p$ and $p'$ relatively prime $\chi_{\Lambda}(x;q)$ is 
a character of a (special) admissible representation 
of A$_{n-1}^{(1)}$, \eqref{KW} corresponding to (a special case of)
the Kac--Wakimoto character formula for admissible characters \cite{KW88}.
The corresponding string functions at fractional-level $N$
are defined via the formal expansion
\begin{equation}\label{Cpp}
\chi_{\Lambda}(x;q)=
\sum_{\mu\in P}q^{\frac{1}{2N}\ip{\mu}{\mu}}C_{\mu,\Lambda}(q) x^{\mu}=
q^{m_{\Lambda}}\sum_{\mu\in P}\C_{\mu,\Lambda}(q) x^{\mu}
\end{equation}
for $\Lambda\in P^{(p,p')}$ with $N$ given by \eqref{Npp}.
To see this is consistent with \eqref{CN}, note that for
$\Lambda\in P_{+}^N$
\begin{multline*}
\sum_{\mu\in P}q^{\frac{1}{2N}\ip{\mu}{\mu}}C_{\mu,\Lambda}(q)x^{\mu}
=\sum_{\mu\in P}\sum_{\substack{\nu\in P/NQ \\ \nu\equiv \mu\pmd{NQ}}} 
q^{\frac{1}{2N}\ip{\mu}{\mu}} C_{\nu,\lambda}(q)x^{\mu} \\
=\sum_{\nu\in P/NQ} C_{\nu,\Lambda}(q)
\sum_{\mu\in NQ+\nu} q^{\frac{1}{2N}\ip{\mu}{\mu}}x^{\mu}
=\sum_{\nu\in P/NQ} C_{\nu,\Lambda}(q)\Theta_{\nu,N}(x),
\end{multline*}
where the first equality follows from \eqref{s1} with $\sigma=(1,2\dots,n)$.
Since the string functions for non-integral levels only satisfy
\eqref{s1} with $\alpha=\emptyset$, it is generally not possible to
rewrite \eqref{Cpp} as \eqref{CN}.

Having defined the configuration sums and string functions for fractional
levels we can now state the following generalization of Conjecture~\ref{conj2}.
\begin{theorem}\label{thm5} 
If Corollary~\ref{corinv} is true then Conjecture~\ref{conj2} holds for all 
$\Lambda\in P^{(p,p')}$. 
\end{theorem}
This theorem claims a conjugate Bailey pair relative to $q^{\ell}$
with $\delta_{\eta}=X_{\eta,\Lambda,N\Lambda_0}(q)$.
It is not clear a priori that substituting this $\delta_{\eta}$ in
\eqref{gddef} will lead to a converging sum, and it is here that 
Conjecture~\ref{conj3} is crucial.

In our proof below we will show that Conjecture~\ref{conj2}
may be manipulated to yield the following expression for
admissible characters of A$_{n-1}^{(1)}$:
\begin{equation}\label{chiX}
\chi_{\Lambda}(q)=q^{m_{\Lambda}}
\sum_{\substack{\mu\in P\\ \mu-\bar{\Lambda}\in Q}} x^{\mu}
\sum_{\substack{\eta\in\Pa\\ \eta_1\leq n-1\\[0.3mm]
|\eta|\bar{\Lambda}_1-\mu\in Q}}
\frac{X_{\eta,\Lambda,N\Lambda_0}(q)}{(q)_{m(\eta)}}
\qbin{m(\eta)}{\frac{|\eta|}{n}(1^n)+\mu}
\end{equation}
for $\Lambda\in P^{(p,p')}$.
Comparing this with \eqref{Cpp} leads to
\begin{equation}\label{CX}
\C_{\mu,\Lambda}(q)=
\sum_{\substack{\eta\in\Pa\\ \eta_1\leq n-1\\[0.3mm]
|\eta|\bar{\Lambda}_1-\mu\in Q}}
\frac{X_{\eta,\Lambda,N\Lambda_0}(q)}{(q)_{m(\eta)}}
\qbin{m(\eta)}{\frac{|\eta|}{n}(1^n)+\mu},
\end{equation}
for $\mu-\bar{\Lambda}\in Q$ and zero otherwise,
from which Theorem~\ref{thm5} easily follows.

\begin{proof}[Proof of Theorem~\ref{thm5}]
Take Corollary~\ref{corinv} and expand the Kostka polynomial therein
using \eqref{Ksup}. By $\tau(\delta)-\delta=\tau(\rho)-\rho$ this yields
\begin{multline*}
\sum_{\substack{\eta\in\Pa \\ \eta_1\leq n-1\\[0.3mm]
|\eta|\equiv |\mu|\pmd{n}}}
\sum_{\sigma,\tau\in S_n}\sum_{\lambda\in nQ+\sigma(\rho)-\rho}
\epsilon(\sigma)\epsilon(\tau)q^{\frac{1}{2n}\ip{\lambda}{\lambda+2\rho}}  \\
\times \frac{1}{(q)_{m(\eta)}}
\qbin{m(\eta)}{\frac{|\eta|-|\nu|}{n}(1^n)+\tau(\nu+\rho)-\rho}
\qbin{m(\eta)}{\frac{|\eta|-|\mu|}{n}(1^n)+\mu-\lambda}
=\delta_{\mu,\nu},
\end{multline*}
where $\mu,\nu\in\Pa$ such that $l(\mu),l(\nu)\leq n$ and $|\mu|=|\nu|$.

Next replace $\mu,\nu$ by $\hat{\mu},\hat{\nu}\in P_{+}$
via $\hat{\omega}=\omega-|\omega|(1^n)/n$, so that
$\hat{\mu}-\hat{\nu}\in Q$.
Making these replacements in the actual equation is trivial and it only 
needs to be remarked that the restriction $|\eta|\equiv|\mu|\pmod{n}$ 
becomes $|\eta|\bar{\Lambda}_1-\hat{\mu}\in Q$. 
After dropping the hats from $\hat{\mu}$ and $\hat{\nu}$ we thus find
\begin{multline}\label{domin}
\sum_{\substack{\eta\in\Pa\\ \eta_1\leq n-1\\[0.3mm]
|\eta|\bar{\Lambda}_1-\mu\in Q}}
\sum_{\sigma,\tau\in S_n}\sum_{\lambda\in nQ+\sigma(\rho)-\rho}
\epsilon(\sigma)\epsilon(\tau)q^{\frac{1}{2n}\ip{\lambda}{\lambda+2\rho}}  \\
\times \frac{1}{(q)_{m(\eta)}}
\qbin{m(\eta)}{\frac{|\eta|}{n}(1^n)+\tau(\nu+\rho)-\rho}
\qbin{m(\eta)}{\frac{|\eta|}{n}(1^n)+\mu-\lambda}
=\delta_{\mu,\nu},
\end{multline}
where $\mu,\nu\in P_{+}$ such that $\mu-\nu\in Q$.
Next we will show that if these conditions are relaxed to
$\mu,\nu\in P$ such that $\mu-\nu\in Q$, then the right-hand side
of \eqref{domin} needs to be replaced by
\begin{equation}\label{rhsnew}
\sum_{\sigma\in S_n} \epsilon(\sigma)\delta_{\mu,\sigma(\nu+\rho)-\rho}.
\end{equation}
In the following we refer to the identity obtained by equating
the left-side of \eqref{domin} with \eqref{rhsnew}
(for $\nu,\mu\in P$ such that $\mu-\nu\in Q$) by (I).
Now there are three possibilities for $\mu,\nu\in P$: 
(i) there exist $w,v\in S_n$
such that $w(\mu+\rho)-\rho,v(\nu+\rho)-\rho\in P_{+}$,
(ii) there exists no such $w$, (iii) there exists no such $v$.
Of course, when $w$ and $v$ exist they are unique and
$w$ $(v)$ does not exist if and only if not all components
of $w+\rho$ ($v+\rho$) are distinct.

First assume case (i) and set $\mu=w^{-1}(\hat{\mu}+\rho)-\rho$,
$\nu=v^{-1}(\hat{\nu}+\rho)-\rho$ with $\hat{\mu},\hat{\nu}\in P_{+}$
and replace $\mu$ and $\nu$ in favour of $\hat{\mu}$ and $\hat{\nu}$.
Then change $\lambda\to w^{-1}(\lambda+\rho)-\rho$,
$\sigma\to w^{-1}\circ\sigma$ and $\tau\to v\circ\tau$ on the left 
and use the symmetry \eqref{Ssymm}. 
On the right change $\sigma\to w^{-1}\circ\sigma\circ v$.
Finally drop the hats of $\hat{\mu}$ and $\hat{\nu}$.
The result of all this is that (I) is transformed into itself
but now $\mu,\nu\in P_{+}$.
For such $\mu$ and $\nu$ the only contributing term 
in the sum \eqref{rhsnew} corresponds to $\sigma=(1,2,\dots,n)$ and 
we are back at \eqref{domin}. Consequently (I) is true in case (i).
Next assume (ii) and replace $\lambda\to \sigma(\lambda+\rho)-\rho$
on the left and $\sigma\to\sigma^{-1}$ on both sides of (I). After
also using \eqref{Ksymm}, the left side then contains the sum
\begin{equation*}
\sum_{\sigma\in S_n}\epsilon(\sigma)
\qbin{m(\eta)}{\frac{|\eta|}{n}(1^n)+\sigma(\mu+\rho)-\rho-\lambda},
\end{equation*}
whereas the right side reads
$\sum_{\sigma\in S_n} \epsilon(\sigma)\delta_{\sigma(\mu+\rho)-\rho,\nu}$.
Since not all components of $\mu+\rho$ are distinct this implies that
both sides vanish. Hence (II) is also true in case (ii).
Similarly it follows that (I) is true in case (iii), with both sides
again vanishing.

To proceed, take identity (I), multiply both sides by $x^{\mu+\rho}$,
sum $\mu$ over $P$ subject to the restriction $\mu-\nu\in Q$ and
finally shift $\nu\to \nu-\rho$.
Adopting definition \eqref{a} for $\nu\in P$ this leads to
\begin{multline*}
\sum_{\substack{\mu\in P\\ \mu-\nu+\rho\in Q}} x^{\mu+\rho}
\sum_{\substack{\eta\in\Pa\\ \eta_1\leq n-1\\[0.3mm]
|\eta|\bar{\Lambda}_1-\mu\in Q}}
\sum_{\sigma,\tau\in S_n}\sum_{\lambda\in nQ+\sigma(\rho)-\rho}
\epsilon(\sigma)\epsilon(\tau)q^{\frac{1}{2n}\ip{\lambda}{\lambda+2\rho}} \\
\times\frac{1}{(q)_{m(\eta)}}
\qbin{m(\eta)}{\frac{|\eta|}{n}(1^n)+\tau(\nu)-\rho}
\qbin{m(\eta)}{\frac{|\eta|}{n}(1^n)+\mu-\lambda}=
a_{\nu}(x),
\end{multline*}
where $\nu\in P$.
We now interchange the sums over $\mu$ and $\lambda$ and replace 
$\mu\to\mu+\lambda$ followed by $\lambda\to \lambda-\rho$.
The left side then factorizes into a double sum over 
$\sigma$ and $\lambda$ and a triple sum over $\mu,\eta$ and $\tau$.
Recalling the theta function \eqref{Tdef}, the double sum
can be recognized as the denominator of the right side of \eqref{KW},
which will be abbreviated here to $A_{\rho}(x;q)$. Hence
\begin{multline*}
q^{-\frac{1}{2n}||\rho||^2}A_{\rho}(x;q)
\sum_{\substack{\mu\in P\\ \mu-\nu+\rho\in Q}} x^{\mu}
\sum_{\substack{\eta\in\Pa\\ \eta_1\leq n-1\\[0.3mm]
|\eta|\bar{\Lambda}_1-\mu\in Q}}
\sum_{\tau\in S_n}\epsilon(\tau) \\
\times \frac{1}{(q)_{m(\eta)}}
\qbin{m(\eta)}{\frac{|\eta|}{n}(1^n)+\tau(\nu)-\rho}
\qbin{m(\eta)}{\frac{|\eta|}{n}(1^n)+\mu}=a_{\nu}(x).
\end{multline*}
Now fix integers $p\geq 1$ and $p'\geq n$,
define $N$ by \eqref{Npp} and let $\Lambda\in P^{(p,p')}$.
Then multiply both sides of the above equation by $q^{\ip{\nu}{\nu}/(2(n+N))}$ 
and sum over $\nu$ such that $\nu\in p'Q+\bar{\Lambda}+\rho$.
On the left swap the sum over $\nu$ with the other three sums and
replace $\nu\to \tau^{-1}(\nu)$.
Recalling the Kac--Wakimoto character formula \eqref{KW} and the definition
\eqref{Xp} of the configuration sum $X_{\eta,\Lambda,\Lambda'}$, this yields 
\eqref{chiX} and thus \eqref{CX}.
To complete the proof choose $\mu=\ell\bar{\Lambda}_1-\pi(k)$, and
use \eqref{Ssymm} with 
$\sigma=\pi$ and $\pi(\ell\bar{\Lambda}_1-\pi(k))=-k-\ell\Lambda_{n-1}$.
Then \eqref{CX} transforms into the claim of the theorem.
\end{proof}

\section{Applications}\label{sec6}
In this final section we give several applications of our A$_{n-1}$
Bailey lemma, resulting in some new A-type $q$-series identities.
Following the derivation of the Rogers--Ramanujan identities outlined
in the introduction, we should try to find pairs of sequences 
$(\alpha,\beta)$ that satisfy the defining relation \eqref{abdef} of a 
Bailey pair.
If these, together with the A$_{n-1}$ conjugate Bailey pairs of
Theorem~\ref{thm1} or \ref{thm5}, are substituted in \eqref{abcdgen}
some hopefully interesting $q$-series identities will result.

According to \eqref{abdef}, what is needed to obtain Bailey pairs are 
polynomial identities involving the antisymmetric supernomials.
Examples of such identities, from which Bailey pairs may indeed be 
extracted, are \eqref{SS}--\eqref{X0}. 
To however explicitly write down the Bailey pairs arising from these 
identities is rather cumbersome due to the fact that we first have
to apply 
\begin{equation*}
\sum_{k\in P} f_k = \sum_{k\in P_{+}}\sum_{\tau\in S_n/S_n^k}f_{\tau(k)}
\end{equation*}
to rewrite the sum over $k$ in \eqref{XN} as a sum over $k\in Q\cap P_{+}$. 
Moreover, once the substitutions in \eqref{abcdgen} have been made, 
one invariably wants to simplify the resulting identity by writing 
the sum over $k\in Q\cap P_{+}$ on the left as a sum over $k\in Q$. 
In order to eliminate the undesirable intermediate step where all 
sums are restricted to $k\in Q\cap P_{+}$, we will never explicitly 
write down Bailey pairs, but only work with the polynomial identities 
that imply these pairs.
To illustrate this in the case of the first Rogers--Ramanujan identity, 
let us take the Bailey pair \eqref{BPR} for $a=1$ and substitute this in
the defining relation \eqref{BP} of an A$_1$ Bailey pair. After a 
renaming of variables this leads to the polynomial identity
\begin{equation}\label{pol}
\sum_{j=-L}^L (-1)^j q^{j(3j-1)/2}\qbin{2L}{L-j}_q=\frac{(q)_{2L}}{(q)_L}.
\end{equation}
Next take the conjugate Bailey pair \eqref{gdr12} for $a=1$ and
substitute this into \eqref{CBP} to find
\begin{equation}\label{jp}
\sum_{r=j}^{\infty} \frac{q^{r^2}}{(q)_{2r}}\qbin{2r}{r-j}_q=
\frac{q^{j^2}}{(q)_{\infty}}
\end{equation}
for $j\in\Z_{+}$. Since both sides are even functions of $j$
we may assume $j\in\Z$.
Now multiply \eqref{jp} by $(-1)^j q^{j(3j-1)/2}$ and sum $j$ over 
the integers. After a change in the order of summation, \eqref{pol} 
may be applied to carry out the $j$-sum on the left, yielding
\begin{equation*}
\sum_{r=0}^{\infty} \frac{q^{r^2}}{(q)_r}=
\frac{1}{(q)_{\infty}}\sum_{j=-\infty}^{\infty}
(-1)^j q^{j(5j-1)/2}.
\end{equation*}
This is exactly the identity obtained when \eqref{BPR} for $a=1$ is
substituted into \eqref{BL}.

In the following we will mimic the above manipulations, albeit in a
more general setting.
First, if we substitute the conjugate Bailey pair of Theorem~\ref{thm5}
in the defining relation \eqref{gddef} and generalize the symmetrization 
immediately after \eqref{jp}, we arrive back at \eqref{CX} which shall 
therefore serve as our starting point. 
Now let $p,p',t,t'$ be integers such that 
$p,t\geq 1$ and $p',t'\geq n$, and set $N=p'/p-n$ and $M=t'/t-n$.
Take \eqref{CX} with $\Lambda\in P^{(t,t')}$,
replace $\mu$ by $k-\rho-\Lambda''$ and multiply both sides by
$\epsilon(\sigma)q^{-\bar{m}_{\Lambda'}
+\ip{k}{k}/(2(n+N))}$ where $\Lambda',\Lambda''\in P^{(p,p')}$ and
$\sigma\in S_n$.
Then sum over $k$ such that $k\in p'Q+\sigma(\bar{\Lambda}'+\rho)$
and sum $\sigma$ over the elements of $S_n$. 
Using the definition \eqref{Xp} of $X_{\eta,\Lambda',\Lambda''}$ this yields
\begin{equation}\label{CXX}
b_{\Lambda,\Lambda',\Lambda''}(q)=
\sum_{\substack{\eta\in\Pa\\ \eta_1\leq n-1\\[0.3mm]
|\eta|\bar{\Lambda}_1-\bar{\Lambda}\in Q}}
\frac{X_{\eta,\Lambda,M\Lambda_0}(q) X_{\eta,\Lambda',\Lambda''}(q)}
{(q)_{m(\eta)}}
\end{equation}
with $b_{\Lambda,\Lambda',\Lambda''}$ defined as
\begin{equation}\label{bll}
b_{\Lambda,\Lambda',\Lambda''}(q)=
q^{-\bar{m}_{\Lambda'}}\sum_{\sigma\in S_n}\epsilon(\sigma)
\sum_{k\in p'Q+\sigma(\bar{\Lambda}'+\rho)}
q^{\frac{\ip{k}{k}}{2(n+N)}}\C_{k-\bar{\Lambda}''-\rho,\Lambda}(q)
\end{equation}
and $\Lambda\in P^{(t,t')}$, $\Lambda',\Lambda''\in P^{(p,p')}$
such that $\bar{\Lambda}-\bar{\Lambda}'+\bar{\Lambda}''\in Q$.

It is interesting to note that \eqref{CXX} implies the highly nontrivial
symmetry 
\begin{equation}\label{bsym1}
b_{\Lambda,\Lambda',N\Lambda_0}=b_{\Lambda',\Lambda,M\Lambda_0}.
\end{equation}
Another symmetry, which follows from \eqref{CXX}, \eqref{Xsymm}
and $(q)_{m(\eta)}=(q)_{m(\eta^c)}$ is 
\begin{equation}\label{bsym2}
b_{\Lambda,\Lambda',\Lambda''}=b_{\Lambda^c,(\Lambda')^c,(\Lambda'')^c}.
\end{equation}
Unlike \eqref{bsym1}, however, this also follows from \eqref{bll}.
Since $\bar{\Lambda}+\rho=-\pi(\bar{\Lambda}^c+\rho)$
it follows from \eqref{KW} that
$\chi_{\Lambda}(x;q)=\chi_{\Lambda^c}(1/x;q)$
and hence that $\C_{\mu,\Lambda}=\C_{-\mu,\Lambda^c}$.
Using this in \eqref{bll}, replacing $k\to -\pi(k)$ and 
$\sigma\to \pi\circ\sigma\circ\pi$ and using 
$\C_{\mu,\Lambda}=\C_{\pi(\mu),\Lambda}$, \eqref{bsym2}
once more follows.

When $\Lambda\in P_{+}^M$, \eqref{bll} coincides
with a special case of an expression of Kac and Wakimoto for the coset 
branching function $b_{\Lambda''}^{\Lambda\otimes\Lambda'}$ defined by
\begin{equation*}
\chi_{\Lambda}(x;q)\chi_{\Lambda'}(x;q)
=\sum_{\substack{\Lambda''\in P^{(p,p'+Mp)} \\[1mm]
\bar{\Lambda}+\bar{\Lambda}'-\bar{\Lambda}''\in Q}}
q^{\phi_{\Lambda,\Lambda',\Lambda''}}
b_{\Lambda''}^{\Lambda\otimes\Lambda'}(q)\chi_{\Lambda''}(x;q)
\end{equation*}
for $\Lambda\in P_{+}^M$, $\Lambda'\in P^{(p,p')}$ and 
\begin{equation*}
\phi_{\Lambda,\Lambda',\Lambda''}=
\frac{||(p'+Mp)(\bar{\Lambda}'+\rho)-p'(\bar{\Lambda}''+\rho)||^2}
{2Mp'(p'+Mp)}
-\frac{||\bar{\Lambda}'-\bar{\Lambda}''||^2}{2M}+m_{\Lambda}.
\end{equation*}
According to \cite[Eq. (3.1.1)]{KW90},
\begin{equation*}
b_{\Lambda''}^{\Lambda\otimes\Lambda'}(q)=
q^{-\bar{m}_{\Lambda'}}\sum_{\sigma\in S_n}\epsilon(\sigma)
\sum_{k\in p'Q+\sigma(\bar{\Lambda}'+\rho)} q^{\frac{\ip{k}{k}}{2(n+N)}}
\C_{\bar{\Lambda}''+\rho-k,\Lambda}(q).
\end{equation*}
Recalling that $\C_{\mu,\Lambda}=\C_{-\mu,\Lambda^c}$ and
comparing the expression for the branching function with \eqref{bll}
leads to the identification
\begin{equation}\label{bb}
b_{\Lambda,\Lambda',\Lambda''}=
b_{M\Lambda_0+\Lambda''}^{\Lambda^c\otimes\Lambda'}
\end{equation}
for $\Lambda\in P_{+}^M$ and $\Lambda',\Lambda''\in P^{(p,p')}$.
Note in particular that $M\Lambda_0+\Lambda''\in P^{(p,p'+Mp)}$
and $\bar{\Lambda}-\bar{\Lambda}'+\bar{\Lambda}''\in Q
\Leftrightarrow \bar{\Lambda}^c+\bar{\Lambda}'-\bar{\Lambda}''\in Q$
as it should.

Thus far we have only manipulated the conjugate Bailey pair of 
Theorem~\ref{thm5}, and to turn \eqref{CXX} into more explicit 
$q$-series identities we need to specify $(p,p')$ and/or $(t,t')$
such that the configuration sums on the right can be evaluated
in closed form. As a first example choose $(p,p')=(1,n)$. 
Then $\Lambda',\Lambda''\in P_{+}^0$ and hence 
$\bar{\Lambda}'=\bar{\Lambda}''=\emptyset$.
Recalling \eqref{X0} and Lemma~\ref{lem2} it follows that the right side 
of \eqref{CXX} trivializes to $\delta_{\Lambda,M\Lambda_0}$.
Making the change $k\to\sigma(nk+\rho)$ in \eqref{bll}
and using $\C_{\mu,\Lambda}=\C_{\sigma(\mu),\Lambda}$ then yields
\begin{equation}\label{AnE}
\sum_{k\in Q}q^{\frac{1}{2}\ip{k}{nk+2\rho}}
\sum_{\sigma\in S_n}\epsilon(\sigma)
\C_{nk+\rho-\sigma(\rho),\Lambda}(q)
=\delta_{\Lambda,M\Lambda_0}
\end{equation}
for $\Lambda\in P^{(t,t')}$ such that $\bar{\Lambda}\in Q$.
This is an A$_{n-1}$, fractional-level analogue of
Euler's pentagonal number theorem.
If for $\Lambda\in P_{+}^1$ we use \eqref{level1}
and perform the sum over $\sigma$ by the Vandermonde determinant 
\eqref{Vdet} in the form
\begin{equation}\label{Vdet2}
\sum_{\sigma\in S_n}\epsilon(\sigma) x^{\sigma(\rho)-\rho}
=\prod_{1\leq i<j\leq n}(1-x_j/x_i).
\end{equation}
with $x_i=q^{-nk_i-n+i}$, we find
\begin{equation*}
\sum_{k\in Q}q^{\frac{1}{2}n(n+1)\ip{k}{k}+\ip{k}{\rho}}
\prod_{1\leq i<j\leq n}\bigl(1-q^{n(k_i-k_j)+j-i}\bigr)=(q)_{\infty}^{n-1}.
\end{equation*}
For $n=2$ this is Euler's famous sum
\begin{equation*}
\sum_{k=-\infty}^{\infty} (-1)^k q^{k(3k+1)/2}=(q)_{\infty},
\end{equation*}
and for general $n$ it corresponds to a specialization of the A$_{n-1}$ 
Macdonald identity \cite{Macdonald72}.
For $\Lambda\in P_{+}^N$ identity
\eqref{AnE} is \cite[Eq. (2.1.17); $\lambda=k\Lambda_0$]{KW90} of 
Kac and Wakimoto, and for $\Lambda\in P^{(t,t')}$ with $n=2$
it is \cite[Prop 8.2; $\eta=0$]{SW00a}.

Many more identities arise if we evoke Theorem~\ref{thm3}.
In fact, instead of just this theorem we will apply the
following generalization which assumes the notation of \eqref{F}.
\begin{conjecture}\label{conj4}
Let $\eta\in\Pa$ such that $\eta_1\leq n-1$.
For $s,u\in\{0,\dots,N\}$, $r,v\in\{0,\dots,n-1\}$ let
$\Lambda=(N-s)\Lambda_r+s\Lambda_{r+1}$ and  
$\Lambda'=(N-u)\Lambda_0+u\Lambda_{n-v}$
such that $|\eta|\bar{\Lambda}_1-\bar{\Lambda}+\bar{\Lambda}'\in Q$
\textup{(}$|\eta|\equiv rN+s+uv\pmod{n}$\textup{)},
and define
\begin{equation*}
F_{\eta,\Lambda,\Lambda'}(q)=
q^{\frac{||\mu-\bar{\Lambda}'||^2-||\bar{\Lambda}||^2}{2N}}
\sum_m q^{\frac{1}{2}m (C\otimes\bar{C}^{-1})m
-(e_1\otimes \bar{e}_{N-s})(I\otimes\bar{C}^{-1})m}\qbin{m+p}{m}_q,
\end{equation*}
where $p$ is determined by
\begin{equation*}
(C\otimes \bar{I})m+(I\otimes\bar{C})p
=(m(\eta)\otimes\bar{e}_1)+(e_1\otimes\bar{e}_{N-s})+(e_v\otimes\bar{e}_{N-u})
\end{equation*}
and the sum is over $m$ such that
\begin{equation*}
\sum_{a=1}^{n-1}\sum_{j=1}^{N-1}jm_j^{(a)}\alpha_a
\in NQ+\mu-(rN+s)\bar{\Lambda}_1-u\bar{\Lambda}_v
\end{equation*}
with $\mu=\sum_{a=1}^{n-1}m_a(\eta)\bar{\Lambda}_a\in P_{+}$. Then 
\begin{equation*}
X_{\eta,\Lambda,\Lambda'}(q)=F_{\eta,\Lambda,\Lambda'}(q).
\end{equation*}
\end{conjecture}
With $\Lambda$ as given above,
$F_{\eta,\Lambda,N\Lambda_0}=F_{\eta,(sN^r)}$ and in this case
the conjecture becomes Theorem~\ref{thm3}.
For $n=2$ Conjecture~\ref{conj4} follows from polynomial identities proved
in \cite{Berkovich94,Warnaar96}.

From \eqref{CXX}, \eqref{bb} and Conjecture~\ref{conj4} we infer the
following fermionic representation for the branching functions:
\begin{equation*}
b^{\Lambda^c\otimes\Lambda'}_{M\Lambda_0+\Lambda''}(q)=
b^{\Lambda\otimes(\Lambda')^c}_{M\Lambda_0+(\Lambda'')^c}(q)=
\sum_{\substack{\eta\in\Pa\\ \eta_1\leq n-1\\[0.3mm]
|\eta|\equiv bM+a\pmd{n}}}
\frac{F_{\eta,\Lambda,M\Lambda_0}(q)F_{\eta,\Lambda',\Lambda''}(q)}
{(q)_{m(\eta)}}
\end{equation*}
with $\Lambda=(M-a)\Lambda_b+a\Lambda_{b+1}$, 
$\Lambda'=(N-s)\Lambda_r+s\Lambda_{r+1}$ and
$\Lambda''=(N-u)\Lambda_0+u\Lambda_{n-v}$ such that
$\bar{\Lambda}-\bar{\Lambda}'+\bar{\Lambda}''\in Q$
(i.e., $bM+a\equiv rN+s+uv\pmod{n}$).
Here $b,r,v\in\{0,\dots,n-1\}$,
$a\in\{0,\dots,M\}$ and $s,u\in\{0,\dots,N\}$.
Note that for $\Lambda''=N\Lambda_0$ this result is certainly true 
thanks to Theorem~\ref{thm3}.
For $b=r=s=0$ (so that $a\equiv uv\pmod{n}$) the above identity has also been 
obtained by Schilling and Shimozono \cite[Eq. (7.10)]{SS01} who exploited 
the fact that
\begin{equation*}
b^{(M-a)\Lambda_0+a\Lambda_1\otimes N\Lambda_0}_{
(M+N-u)\Lambda_0+u\Lambda_v}(q)=
\lim_{L\to\infty}q^{-akL-nM\binom{kL}{2}}
K^{N+M}_{\lambda^{(L)},R^{(L)}}(q).
\end{equation*}
Here $K^N_{\lambda,R}(q)$ is the level-restricted generalized
Kostka polynomial, $k\in [n-1]$, $\lambda^{(L)}=(c^{n-v}(c+u)^v)$
with $c=kLM+(a-uv)/n$, and $R^{(L)}=(R_1,\dots,R_{Ln+1})$,
with $R_1=(a)$ and $R_i=(M^k)$ for $i\geq 2$.
By using the generalization of Theorem~\ref{thm3} to all 
$\Lambda\in P_{+}^N$ \cite{OSS01,SS00,SS01} we can also
apply \eqref{CXX} to prove fermionic formulae for all of the branching 
functions $b^{\Lambda\otimes\Lambda'}_{(N+M)\Lambda_0}$ with 
$\Lambda\in P_{+}^M$ and $\Lambda'\in P_{+}^N$.

Another nice specialization of \eqref{CXX} arises for $\Lambda\in P_{+}^1$
and $\Lambda''=N\Lambda_0$.
Then $\bar{\Lambda}-\bar{\Lambda}'\in Q$ fixes $\Lambda$, and
by \eqref{X1}, \eqref{level1} and \eqref{Vdet2} it follows that
all dependence on $\Lambda$ cancels out.
Dropping the prime in $\Lambda'$, this yields
\begin{multline}\label{prodX}
\frac{q^{\frac{1}{2}||\bar{\Lambda}||^2}}
{(q)_{\infty}^{n-1}}\sum_{k\in Q}
q^{\frac{1}{2}p'(p+p')\ip{k}{k}+(p+p')\ip{k}{\bar{\Lambda}}+p\ip{k}{\rho}}\\
\times \prod_{1\leq i<j\leq n}
\bigl(1-q^{p'(k_i-k_j)+j-i+\ip{\bar{\Lambda}}{\varepsilon_i-\varepsilon_j}}
\bigr)=\sum_{\substack{\eta\in\Pa\\ \eta_1\leq n-1\\[0.3mm]
|\eta|\bar{\Lambda}_1-\bar{\Lambda}\in Q}} \!\! 
\frac{q^{\frac{1}{2}\ip{\mu}{\mu}}
X_{\eta,\Lambda,N\Lambda_0}(q)}{(q)_{m(\eta)}}
\end{multline}
for $\Lambda\in P^{(p,p')}$ and 
$\mu=\sum_{a=1}^{n-1}m_a(\eta)\bar{\Lambda}_a$.
When $p'=n$ the sum on the left can be performed by the
A$_{n-1}$ Macdonald identity, resulting in 
\begin{equation}\label{Xnp}
\sum_{\substack{\eta\in\Pa\\ \eta_1\leq n-1\\[0.3mm]
|\eta|\equiv 0\pmd{n}}}
\frac{q^{\frac{1}{2}\ip{\mu}{\mu}}X_{\eta}^{p,n}(q)}
{(q)_{m(\eta)}}
=\frac{(q^{n+p};q^{n+p})_{\infty}^{n-1}}{(q)_{\infty}^{n-1}}
\prod_{i=1}^{n-1}(q^i,q^{n+p-i};q^{n+p})_{\infty}^{n-i}
\end{equation}
with 
\begin{equation*}
X_{\eta}^{p,n}(q)=
\sum_{\sigma\in S_n}\epsilon(\sigma) \sum_{k\in nQ}
q^{\frac{p}{2n}\ip{k}{k+2\rho}}\qbin{m(\eta)}
{\frac{|\eta|}{n}(1^n)+k+\rho-\sigma(\rho)}.
\end{equation*}
Here $X_{\eta}^{p,n}(q)=X_{\eta,N\Lambda_0,N\Lambda_0}(q)$ for $N=n(1/p-1)$
and $N\Lambda_0\in P^{(p,n)}$.

The identity \eqref{Xnp} may be viewed as an A$_{n-1}$ version of the
second Rogers--Ramanujan identity and its higher moduli 
generalizations \cite{Andrews74,Bressoud80}.
By standard Bailey chain type arguments 
\cite{Andrews84,Andrews01,W01} it follows that for 
$k\geq 1$ and $a\in\{0,1\}$, 
\begin{equation*}
X^{2k+a,2}_{(1^{2L})}(q)
=\sum_{L\geq n_1\geq \dots\geq n_{k-1}\geq 0}
\frac{q^{L+n_1(n_1+1)+\cdots+n_{k-1}(n_{k-1}+1)}(q)_{2L}}
{(q)_{L-n_1}\cdots(q)_{n_{k-2}-n_{k-1}}(q^{2-a};q^{2-a})_{n_{k-1}}}
\end{equation*}
which substituted into \eqref{Xnp} gives
\begin{multline}\label{AB}
\sum_{n_1\geq \dots\geq n_k\geq 0}
\frac{q^{n_1(n_1+1)+\cdots+n_k(n_k+1)}}
{(q)_{n_1-n_2}\cdots(q)_{n_{k-1}-n_k}(q^{2-a};q^{2-a})_{n_k}}\\
=\frac{(q,q^{2k+a+1},q^{2k+a+2};q^{2k+a+2})_{\infty}}{(q)_{\infty}}.
\end{multline}
There is in fact an alternative route from \eqref{Xnp} to this result,
based on level-rank duality. First observe that if we write the right-hand
side of \eqref{Xnp} as $f_{n,p}(q)$, then $f_{n,p}=f_{p,n}$.
The identities of \eqref{AB} should therefore also result from the
$p=2$ case of \eqref{Xnp}. Indeed, experiments for $n=3$ and $4$
confirm the following conjecture:
\begin{equation*}
X_{\eta}^{2,2k+a}(q)=\frac{q^{n_1+\cdots+n_k}(q)_{m(\eta)}}
{(q)_{n_1-n_2}\cdots(q)_{n_{k-1}-n_k}(q^{2-a};q^{2-a})_{n_k}}
\end{equation*}
if $m(\eta)=(n_1-n_2,\dots,n_{k-1}-n_k,n_k,n_k,n_{k-1}-n_k,\dots,n_1-n_2)$
and zero otherwise. Here $n_k,n_k$ should read $2n_k$ when $a=0$.
Substituted into \eqref{Xnp} this once again yields \eqref{AB}.

As a final application we employ the following conjecture
for $n=3$ and $\Lambda\in P^{(2,4)}_{+}$:
\begin{equation*}
X_{\eta,\Lambda,-\Lambda_0}(q)=q^{\frac{1}{6}(\zeta_1-\zeta_2)^2+
\frac{1}{2}(\zeta_1+\zeta_2)-||\bar{\Lambda}||^2}=
q^{\binom{\zeta_1+1}{2}+\binom{\zeta_2+1}{2}-\frac{1}{2}\zeta C^{-1}\zeta
-||\bar{\Lambda}||^2},
\end{equation*}
where $\zeta=m(\eta)\in\Z_{+}^2$. Since the elements of
$P^{(2,4)}_{+}$ are given by $\Lambda_l-2\Lambda_0$ for $l\in\{0,1,2\}$,
the condition $|\eta|\bar{\Lambda}_1-\bar{\Lambda}\in Q$ is equivalent to
$\zeta_1+2\zeta_2+l\equiv 0\pmod{3}$. 
Substituting this last conjecture into \eqref{prodX} gives
\begin{multline*}
\frac{1}{(q)_{\infty}^2}\sum_{k\in Q}
q^{12\ip{k}{k}+6\ip{k}{\bar{\Lambda}_l}+2\ip{k}{\rho}}
\prod_{1\leq i<j\leq 3}
\bigl(1-q^{4(k_i-k_j)+j-i+\ip{\bar{\Lambda}_l}{\varepsilon_i-\varepsilon_j}}
\bigr) \\
=\sum_{\substack{\zeta_1,\zeta_2\geq 0\\ \zeta_1+2\zeta_2\equiv -l\pmd{3}}}
\frac{q^{\binom{\zeta_1+1}{2}+\binom{\zeta_2+1}{2}
-\frac{3}{2}||\bar{\Lambda}_l||^2}}
{(q)_{\zeta_1}(q)_{\zeta_2}}
\end{multline*}
whereas substitution into \eqref{CXX} leads to
\begin{equation*}
b_{\Lambda_l-2\Lambda_0,\Lambda_l-2\Lambda_0,-2\Lambda_0}(q)=
\sum_{\substack{\zeta_1,\zeta_2\geq 0\\ \zeta_1+2\zeta_2\equiv -l\pmd{3}}}
\frac{q^{\frac{1}{3}(\zeta_1-\zeta_2)^2+\zeta_1+\zeta_2
-2||\bar{\Lambda}_l||^2}}
{(q)_{\zeta_1}(q)_{\zeta_2}}.
\end{equation*}

\subsection*{Acknowledgement}
I thank Will Orrick, Anne Schilling and Trevor Welsh 
for many helpful discussions.

\end{document}